\documentclass[graybox]{svmult}

\usepackage{amsmath, amssymb}
\usepackage{mathptmx}       
\usepackage{helvet}         
\usepackage{courier}        
\usepackage{type1cm}        
\usepackage{makeidx}         
\usepackage{graphicx}        
\usepackage{multicol}        
\usepackage[bottom]{footmisc}

\makeindex             
\newcommand{\bbR}{{\mathbb{R}}}

\newcommand{\nablas}{\nabla_\Gamma}
\newcommand{\nablash}{\nabla_{\Gamma_h}}

\newcommand{\Omegaoh}{\Omega_{h,1}}

\newcommand{\bfn}{\boldsymbol n}
\newcommand{\bfbeta}{\boldsymbol \beta}

\newcommand{\bfx}{\boldsymbol x}
\newcommand{\bfa}{\boldsymbol a}
\newcommand{\bfb}{\boldsymbol b}

\newcommand{\mcK}{\mathcal{K}}

\newcommand{\mcF}{\mathcal{F}}
\newcommand{\mcN}{\mathcal{N}}

\newcommand{\Gammah}{{\Gamma_h}}

\newtheorem{rem}{Remark}[section]

\newcommand{\divs}{\text{div}_\Gamma}
\newcommand{\divsh}{\text{div}_{\Gamma_h}}

\newcommand{\IR}{\mathbb{R}}
 
\usepackage{booktabs}

\begin{document}

\title*{
A Space-Time Cut Finite Element Method with quadrature in time} 

\titlerunning{A Space-Time CutFEM with quadrature in time}
\author{Sara Zahedi}
\institute{Sara Zahedi \at Department of Mathematics, KTH Royal Institute of Technology, SE-100 44 Stockholm, Sweden, \email{sara.zahedi@math.kth.se}}
%
%
\maketitle

\abstract{We consider convection-diffusion problems in time-dependent domains and present a space-time finite element method based on quadrature in time which is simple to implement and avoids remeshing procedures as the domain is moving. 
The evolving domain is embedded in a domain with fixed mesh and a cut finite element method with continuous elements in space and discontinuous elements in time is proposed.  The method allows the evolving geometry to cut through the fixed background mesh arbitrarily and thus avoids remeshing procedures. However, the arbitrary cuts may lead to ill-conditioned algebraic systems. A stabilization term is added to the weak form which guarantees well-conditioned linear systems independently of the position of the geometry relative to the fixed mesh and in addition makes it possible to use quadrature rules in time to approximate the space-time integrals. 
We review here the space-time cut finite element method presented in~\cite{HLZ16ST} where linear elements are used in both space and time and extend the method to higher order elements for problems on evolving surfaces (or interfaces). We present a new stabilization term which also when higher order elements are used controls the condition number of the linear systems from cut finite element methods on evolving surfaces. The new stabilization combines the consistent ghost penalty stabilization~\cite{Bu10} with a term controlling normal derivatives at the interface.}

\section{Introduction}\label{sec:intro}
Finite Element Methods (FEM) are well known for efficiently solving Partial Differential Equations (PDEs) in complex geometries. However, when the geometry is moving a remeshing procedure is needed to fit the mesh to the evolving geometry. 
In for example simulations of multiphase flow phenomena the evolving geometry can be the interface separating two immiscible fluids or the domain occupied by one of the fluids. Topological changes such as drop-breakup or coalescence occur and the remeshing process is both complicated and expensive, especially in three space dimensions. In~\cite{HLZ16ST} and~\cite{HaLaZa15} we therefore present cut finite element methods that, contrary to standard FEM, allow the evolving geometry to be arbitrarily located with respect to a fixed background mesh.

In Cut Finite Element Methods (CutFEM) the domain where the PDE has to be solved is embedded in a computational domain with fixed background mesh equipped with a standard finite element space and one uses the restriction of the basis functions to the so called active mesh where the bilinear forms associated with the weak formulation are evaluated. A stabilization term is added in the weak form to ensure well-conditioned linear systems independently of the position of the geometry relative to the background mesh.  
 
In~\cite{HaLaZa15} the strategy is to follow characteristics to fetch information from interfaces at previous time steps. Error estimates are derived in the L$^2$-norm for a convection diffusion equation on a moving interface. The method is only first order accurate in the L$^2$-norm. In~\cite{HLZ16ST} the strategy is instead to use a space-time finite element method and the present contribution is built on this idea. 
Compared to prior work on Eulerian space-time finite element methods such as~\cite{Gr14, CL15, OlRe14, OlReXu14a} a stabilization term is added to the weak formulation. Due to this stabilization the method in~\cite{HLZ16ST} has the following characteristics: 1) the linear systems resulting from the method have bounded condition numbers independently of how the geometry cuts through the background mesh; 2) the implementation of the method can be based on directly approximating the space-time integrals using quadrature rules for the integrals over time.
Due to the second point, provided that a method for the representation and evolution of the geometry is available, it is straightforward to implement the space-time CutFEM in~\cite{HLZ16ST} starting from a stationary CutFEM. This makes the implementation convenient when going to higher order elements and coupled bulk-surface problems. A space-time unfitted finite element method using the trapezoidal rule to approximate the integral over time was proposed and studied in~\cite{Gr14} but the method failed to converge in case of moving interfaces. We note that no stabilization was used in~\cite{Gr14}.

In this contribution we review the space-time method in~\cite{HLZ16ST} for solving convection-diffusion equations modeling the evolution of surfactants and extend the method to higher order elements for problems on moving interfaces. A new stabilization term is proposed which in contrast to the stabilization term in~\cite{HLZ16ST} leads to linear systems with condition numbers scaling as $\mathcal{O}(h^{-2})$ for linear as well as higher order elements.

The remainder of this contribution is outlined as follows. We start with a surface problem in Section~\ref{sec:model}. We state the surface PDE in~\ref{sec:MathmodS} and present the space-time CutFEM and the new stabilization term in Section~\ref{sec:methodS}. Implementation aspects are discussed in Section~\ref{sec:implementS} and we show numerical examples using both linear and higher order elements in Section~\ref{sec:numexpS}. Next we consider a coupled bulk-surface problem. We present the computational method, implementation aspects, and a numerical example from~\cite{HLZ16ST} in Section~\ref{sec:method}-\ref{sec:numexp}.  We discuss our results in Section~\ref{sec:discussion}.

\section{A Surface Problem}
\label{sec:model}
\begin{figure}
\begin{center}
\includegraphics*[width=0.4\textwidth]{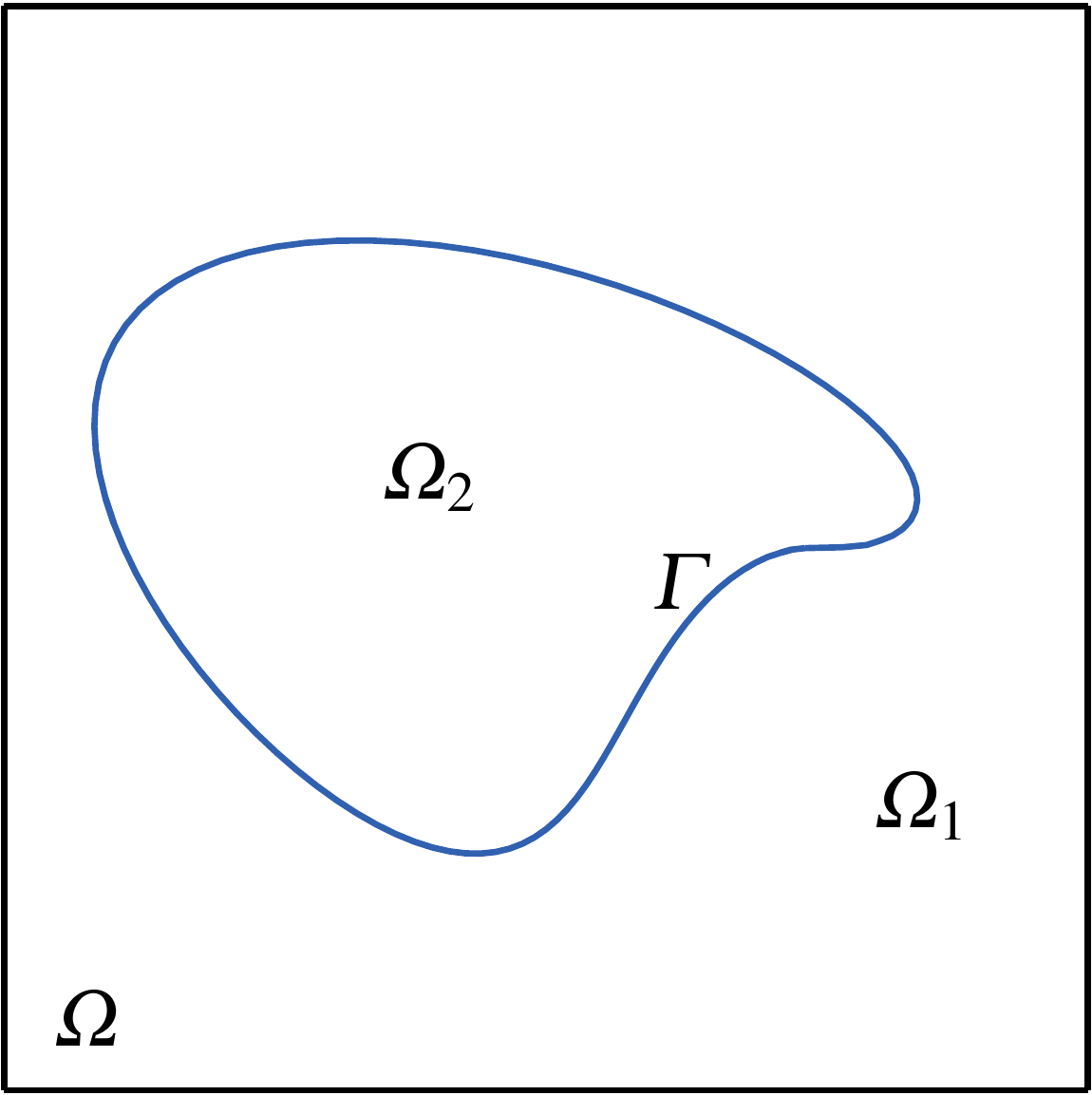} 
\caption{Illustration of the domain $\Omega \in {\bbR}^2$ and the two subdomains $\Omega_i(t)$, $i=1,2$ that are separated by an interface $\Gamma$.\label{fig:domain}}
\end{center}
\end{figure}

Consider an open bounded domain $\Omega$ in ${\bbR}^d$, $d=2,3$ with
convex polygonal boundary $\partial \Omega$. During all time $t$ in
the interval $I=[0,T]$ this domain contains two subdomains
$\Omega_1(t)$ and $\Omega_2(t)$ that are separated by a smooth
interface $\Gamma(t) = \partial \Omega_1(t) \cap \partial
\Omega_2(t)$, a simply connected closed curve in $\IR^2$ or a surface
in $\IR^3$ with exterior unit normal $n\in \bbR^d$. The
interface is moving with a given velocity field $\bfbeta: I \times
\Omega \rightarrow {\bbR}^d$ and does not intersect the boundary of
the domain $\partial \Omega$ ($\Gamma \cap \partial \Omega =
\emptyset$) or itself for any $t\in I$. See Fig.~\ref{fig:domain} for
an illustration in two dimensions.

For $x \in \bbR^d$ let $p(x)$ be the closest point projection mapping onto $\Gamma$. Let $U_{\delta_0}(\Gamma)$ denote the tubular neighborhood of the interface $\Gamma$ in which for each $x\in U_{\delta_0}(\Gamma)$ there is a unique $p(x)$ on $\Gamma$. We may extend any function defined on $\Gamma$ to $U_{\delta_0}(\Gamma)$ by $u^e=u(p(x))$,  $x\in U_{\delta_0}(\Gamma)$. 
We use this extension to for example define the tangential derivative on $\Gamma$ as: 
\begin{equation}
\nabla_\Gamma u=P_\Gamma \nabla u^e 
\end{equation}
where 
\begin{equation}\label{eq:Proj}
P_\Gamma=I-n \otimes n 
\end{equation}
Here $I$ is the identity matrix, and $\otimes$ denotes the outer product $(\bfa \otimes \bfb)_{ij} = a_i b_j$ for any two vectors $\bfa$ and $\bfb$. Note that the tangential derivative depends only on the values of u on $\Gamma$ and does not depend on the particular choice of extension. In the following we will leave the superscript off and write $u$ also for the extended function.

\subsection{Mathematical model}\label{sec:MathmodS}
Consider the following time dependent convection-diffusion equation: 
\begin{equation}
\partial_t u + \bfbeta \cdot \nabla u +(\divs \bfbeta) u  - k_S \Delta_\Gamma u  =  f
\quad \text{on  $\Gamma(t)$, $\quad t \in I$} \label{eq:uS}
\end{equation}
with initial condition $u(0,\bfx) = u^0$ on $\Gamma(0)$ ($x \in \bbR^d$). Here $ f \in L^2(t \times \Gamma(t))$ for all $t\in I$, $k_S>0$ is the diffusion coefficient, $\Delta_\Gamma$ is the Laplace-Beltrami operator, $\partial_t=\frac{\partial}{\partial t}$, and 
\begin{equation}
\divs = \text{tr}((I-n \otimes n)\nabla)
\end{equation}

\begin{rem}
When $f=0$ equation~\eqref{eq:uS} models the evolution of the concentration of insoluble surfactants on an interface separating two immiscible fluids. 
Then the following conservation law also holds:
\begin{equation}\label{eq:conservS} 
\int_{\Gamma(t)} u ds= \overline{u}_0 \quad \textrm{for all $t\geq 0$}
\end{equation} 
\end{rem}

\subsection{The space-time cut finite element method}\label{sec:methodS}
We now propose a space-time cut finite element method for solving the surface PDE stated in the previous section. The method uses the strategy in~\cite{HLZ16ST}.

\subsubsection{Mesh and space}\label{sec:meshspaceSurf}
Create a quasiuniform partition of $\Omega$ into shape regular triangles for $d=2$ and tetrahedra for $d=3$ of diameter $h$ and denote it by  $\mcK_{0,h}$. We will refer to this partition as the fixed background mesh. Let $V_{0,h}^p$ be the space of continuous piecewise polynomials of degree $p \geq 1$ defined on the background mesh $\mcK_{0,h}$. 
Partition the time interval $I=[0, T]$,  $0 = t_0 < t_1 <\dots < t_N = T$, into time steps $I_n = (t_{n-1},t_n]$ of length $k_n = t_n - t_{n-1}$ for $n = 1,2,\dots,N$. 

\begin{figure}\centering
\includegraphics[scale=0.4]{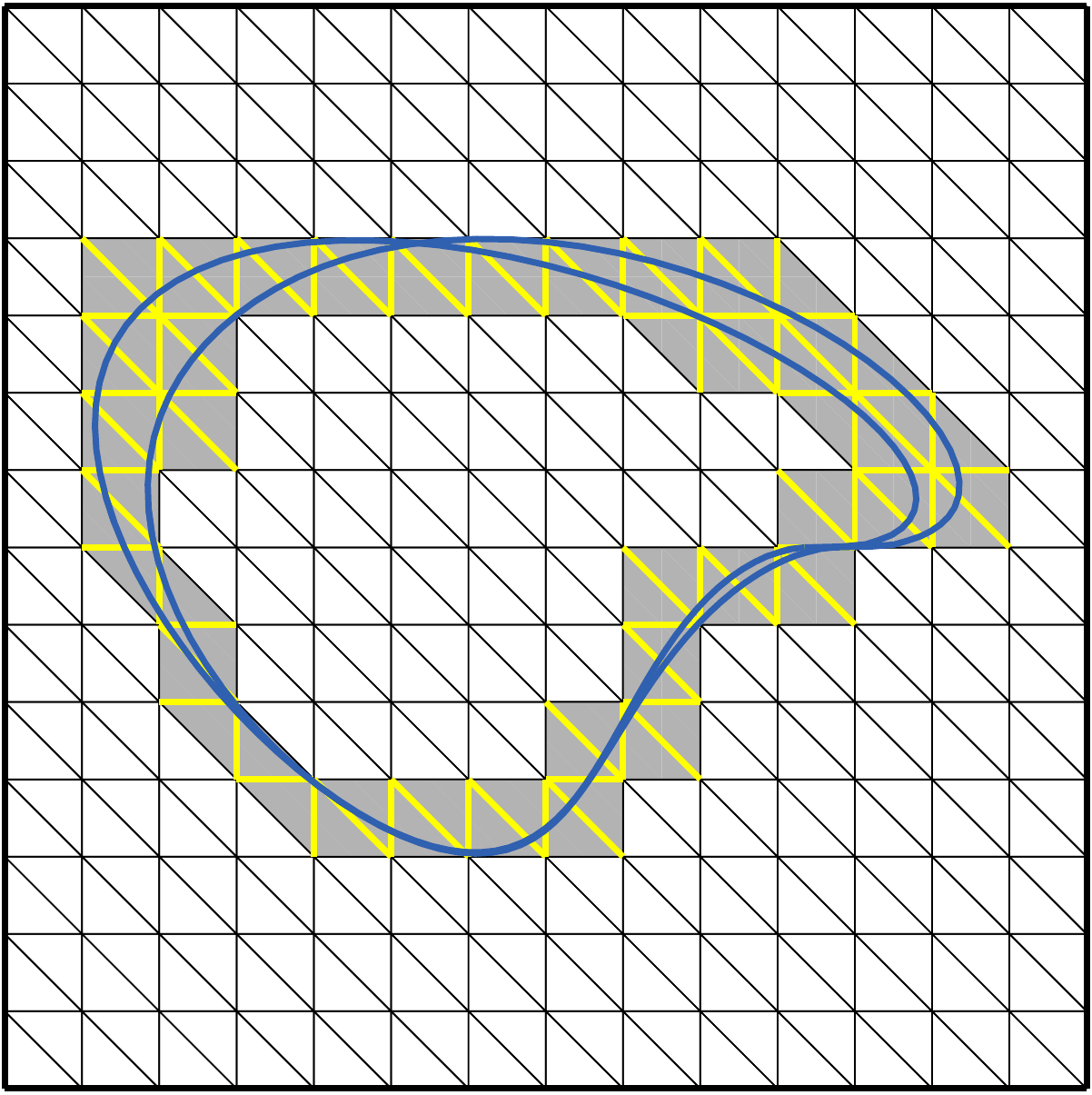} 
\caption{Illustration of the sets introduced in Section~\ref{sec:meshspaceSurf}. The two blue curves show the position of the interface at the endpoints $t=t_{n-1}$ and $t=t_{n}$ of the time interval $I_n = (t_{n-1},t_n]$. The shaded domain shows $\mcN_{S,h}^n$, the so called active mesh. Edges in $\mcF_{S,h}^n$  are marked with yellow thick lines.\label{fig:illustSurf}}.
\end{figure}

Denote the set of elements in the fixed background mesh that are cut by the interface by $\mcK_{S,h}$:
\begin{equation}
\mcK_{S,h}(t) = \{K \in \mcK_{0,h} \,:\, K \cap \Gamma(t) \neq \emptyset \}
\end{equation}
and define the following domain: 
\begin{equation}\label{eq:mcnmcs}
\mcN_{S,h}^n = \bigcup_{t \in I_n} \bigcup_{K \in \mcK_{S,h}(t)} K  
\end{equation}
We refer to the last domain as the active mesh. For an illustration of the active mesh in two dimensions see the shaded domain in Fig.~\ref{fig:illustSurf}.

Associated to the active mesh is the space-time slab $S_S^n = I_n \times \mcN_{S,h}^n$ on which we define the space $V_{S,h}^n$:
\begin{equation}
V_{S,h}^n=P_q(I_n) \otimes V_{0,h}^p|_{\mcN_{S,h}^n}
\end{equation}
Here $P_{q}(I_n)$ is the space of polynomials of degree less or equal to $q$ on $I_n$. 
Functions $v_h(t,x)$ in $V_{S,h}^n$ take the form
\begin{equation}\label{eq:ansatzSurf}
v_h(t,x)=\sum_{j=0}^{q} v_{S,j}(x)\left(\frac{t-t_{n-1}}{k_n} \right)^j 
\end{equation}
where $t\in I_n$ and $v_{S,j}(x)$, $j=0,1,\cdots, q$ are functions in $V_{0,h}^p|_{\mcN_{S,h}^n}$  (the space of restrictions to the active mesh of functions in $V_{0,h}^p$) and hence can be written as 
\begin{align}\label{eq:funcwhSurf}
v_{S,j}(x)= \sum_{i}  \xi_{ij}^S \varphi_i(\bfx)|_{\mcN_{S,h}^n} 
\end{align}
Here $\xi_{ij}^S \in {\bbR}$ are coefficients and $\varphi_i(\bfx)$ are the standard basis functions in the space $V_{0,h}^p$ associated with the degree of freedom $i$. The sum is over all degrees of freedom in the active mesh (the shaded domain in Fig.~\ref{fig:illustSurf}).

\subsubsection{The variational formulation}\label{sec:weakformSurf}
For $t \in I_n$ and given $u_{h}(t_{n-1}^-,x)$ the weak formulation is to find $u_{h} \in V_{S,h}^n$ such that 
\begin{align}\label{eq:spacetimeformS}
A_h^n(u_h,v_h) + J_h^n (u_h,v_h) &=\int_{I_n} (f,v_h)_{\Gamma(t)}
\quad \forall v_h \in V_{S,h}^n
\end{align}
Here 
\begin{equation}
A_h^n (u,v) =  \int_{I_n} (\partial_t u,v)_{\Gamma(t)} \, dt + \int_{I_n} a_h(t,u,v)\, dt 
+([u],v(t_{n-1}^+,\bfx))_{\Gamma(t_{n-1})} 
\end{equation}
with 
\begin{equation} \label{eq:ahsurf}
a_{h}(t,u,v)= (\bfbeta \cdot \nabla u,v)_{\Gamma(t)} +((\divs \bfbeta) u,v)_{\Gamma(t)}  + \left(k_S \nablas u,\nablas v\right)_{\Gamma(t)}
\end{equation}
and $J_h^n$ is a stabilization term we will introduce and discuss in
the next section.  

Note that the trial and test functions are discontinuous from one space-time slab to another and therefore at a given time $t_n$ (where $n$ is the time step number) there are two distinct solutions, at times $t_n^\pm := \lim_{\epsilon\rightarrow 0}t_n\pm \epsilon$. To weakly enforce continuity at $t_n$ the term $([u],v(t_{n-1}^+,\bfx))_{\Gamma(t_{n-1})}$ is added and the discrete equations can be solved one space-time slab at a time, see \emph{e.g.}~\cite{Ja78}.

\subsubsection{Stabilization}
We have two aims with the stabilization term added to the weak form: 1) to control the condition number of the resulting linear systems independently of how the geometry cuts through the background mesh; 2) to be able to directly approximate the space-time integrals using quadrature rules for the integrals over time. We propose the stabilization
 \begin{equation} \label{eq:stabS}
J_h^n(u_h,v_h)= \int_{I_n} j_{h, F}(u_h,v_h)\, dt+ \int_{I_n}  j_{h,\Gamma}(t,u_h,v_h)\, dt
\end{equation}
where we combine the face stabilization
\begin{align}\label{eq:stab1}
j_{h, F}(u_h,v_h) &= \sum_{F \in \mcF_{S,h}^n} \sum_{i=1}^{p} c_{F,i} h^{\gamma}([\partial_n^i u_h]\vert_F,[\partial_n^i v_h]\vert_F)_F
\end{align}
with the interface stabilization
\begin{align}\label{eq:stab2}
j_{h,\Gamma}(t,u_h,v_h)= \sum_{i=1}^{p}c_{\Gamma,i} h^{\gamma}(\partial_n^i u_h,\partial_n^iv_h)_{\Gamma(t)}
\end{align}
Here $\partial_n^i$ denotes the normal derivative of order i, $[x]\vert_F$ denotes the jump of $x$ over the face $F$, $\mcF_{S,h}$ is the set of internal faces, i.e. faces with two neighbors in the active mesh $\mcN_{S,h}^n$, see the yellow marked edges in Fig.~\ref{fig:illustSurf}, $c_{F,i}>0$, $c_{\Gamma,i}>0$ are stabilization constants, and  we take
\begin{equation}
\gamma=2i
\end{equation}
This choice of $\gamma$ yields the weakest stabilization which still controls the condition number. 
 
The ghost penalty stabilization, here referred to as the face stabilization, has been used in several works on CutFEM~\cite{BuHa12, BuHaLa15, BuHaLaZa16, HLZ16ST, HaLaZa14, HaLaZa15}, though originally proposed in~\cite{Bu10}, to control the condition number of the resulting system matrix independently of how elements in the fixed background mesh are cut by the geometry. 

For surface PDEs, adding a face stabilization to the weak form was first proposed in~\cite{BuHaLa15} as a way to get $\mathcal{O}(h^{-2})$ condition number estimates starting from the TraceFEM in~\cite{OlReGr09} for solving the Laplace-Beltrami equation on a stationary surface. Note that in the stabilized method~\cite{BuHaLa15} the active mesh and the finite element spaces  differ from the unstabilized method~\cite{OlReGr09}.  In the TraceFEM in~\cite{OlReGr09} the active mesh is the restriction of the background mesh to the surface which results in an induced cut surface mesh while in the CutFEM in~\cite{BuHaLa15} the active mesh is the union of all elements that are cut by the surface. The finite element spaces are then defined as the restriction of the finite element space defined on the background mesh to the active mesh. The same face stabilization term also leads to a stable discretization for convection-diffusion equations in the case of dominating convection~\cite{BuHaLaZa15} and no other stabilization term such as for example a SUPG term is needed. 

For PDEs on evolving surfaces, to the author's best knowledge, the face stabilization term has only been used with linear elements in space. For higher order elements, following the same scaling as in previous work, the parameter $\gamma$ in the face stabilization, equation~\eqref{eq:stab1}, should be $2i-2$. We will see in Section~\ref{sec:numexpS} that for higher order elements this face stabilization alone does not provide enough control and the condition numbers of the resulting linear systems do not scale as $\mathcal{O}(h^{-2})$. 

We also note that in~\cite{DeElRa14} for linear elements the full gradient on the interface rather than the tangential gradient in equation~\eqref{eq:ahsurf} was used to get control over the normal derivative of the finite element solution. This corresponds to choosing, $p=1$, $\gamma=0$, $c_{\Gamma}^1=1$ and $c_F^1=0$ in~\eqref{eq:stab2}. However, in~\cite{Re15} an example is given that shows that such a surface stabilization does not give $\mathcal{O}(h^{-2})$ condition number estimates  for higher order elements, at least not independently of how the surface cuts the background mesh. 

We propose to combine the face stabilization \eqref{eq:stab1} and the stabilization involving the normal derivatives at the interface \eqref{eq:stab2} and to take $\gamma=2i$. Note that our choice of $\gamma$ gives a different scaling of the face stabilization than what is used in~\cite{BuHaLaZa15} and a different scaling of the interface or surface stabilization term than what is used in~\cite{DeElRa14} for linear elements. The idea with the new stabilization is to use the face stabilization to reach elements which have a large intersection with the interface and on those elements use the interface stabilization term to get enough control. In~\cite{LaZa16} we propose a CutFEM for the Laplace-Beltrami equation on a stationary surface with such a stabilization term and prove that the condition number of the resulting linear system also for higher order elements scales as $\mathcal{O}(h^{-2})$ independently of how the surface cuts the background mesh. Recently, another stabilization, a normal gradient stabilization which acts on the elements in the active mesh has been proposed in~\cite{BuHaLaMa16, GrLeRe17}. For the Laplace-Beltrami equation on a stationary surface this stabilization has also been proven to yield $\mathcal{O}(h^{-2})$ condition number estimates independent of the degree of the polynomials used in the trial and the test space~\cite{GrLeRe17}.

\subsection{Implementation}~\label{sec:implementS}
Often the exact interface $\Gamma$ is not available but an approximation $\Gamma_h$ is.  This means
that in the definition of the active mesh, the finite element spaces, and in all integrals in the weak formulation in 
Section~\ref{sec:methodS} the exact interface $\Gamma$ and the normal $n$ are replaced by an approximate interface $\Gammah$ and normal $n_h$. See Section~\ref{sec:Intrep} on how the exact interface $\Gamma$ is approximated in this work.

We approximate all space-time integrals in the weak form by using first a quadrature rule in time and then a quadrature rule in space. The proposed space-time formulations in~\cite{Gr14, CL15, OlRe14, OlReXu14a} instead convert the space-time integrals to surface integrals over the space-time manifold
\begin{equation}
S=\cup_{t\in (0, T)} \Gamma(t) \times {t}
\end{equation}
by using the identity
\begin{equation}
\int_0^T\int_{\Gamma(t)} f(s,t) \,dsdt=\int_S f(s)\frac{1}{(1+(\beta\cdot \bfn)^2)^{1/2}} \,ds 
\end{equation}
Hence when $\Gamma$ is a surface in ${\bbR}^d$ surface integrals in ${\bbR}^{d+1}$ need to be computed. The space-time manifold $S \subset \bbR^{d+1}$ is approximated by a discrete surface $S_h$ and integrals are computed over $S_h$. 
We propose to directly approximate the space-time integrals using quadrature rules for the integrals over time, see Section~\ref{sec:assembly} for more details. Geometric computations, involving the construction of the interface $\Gamma$, are then done only at the quadrature points in time. This essentially means that it is straightforward to implement the proposed space-time CutFEM starting from a stationary CutFEM. This is possible due to the stabilization term we add to the weak form. 
We note that an advantage of the space-time method in~\cite{OlReXu14a} is the existing analysis~\cite{OlRe14} of the method. However, optimal estimates were proved in a weaker norm than the L$^2$-norm.

\subsubsection{Numerical representation of the interface}~\label{sec:Intrep} 
The interface is represented using either an \emph{explicit}\/ representation, by marker particles and a parametrization, see e.g.~\cite{Pe77}, or an \emph{implicit}\/ representation by the level set of a higher dimensional function, see e.g.~\cite{OsSe88}. In this work we only use the level set method when linear elements in space and time are used,  i.e. $p=q=1$.

\textbf{An implicit representation}
Let the level set function $\phi(t,\bfx)$ be the signed distance function with positive sign in $\Omega_2$. The interface is then defined implicitly as the zero countour of $\phi$. The spatial gradient of the signed distance function defines the exterior unit normal on $\Gamma$ with respect to $\Omega_1$:
\begin{equation}\label{eq:normal}
\bfn(t,\bfx ) = \nabla \phi(t,\bfx )=\frac{\nabla \phi(t,\bfx )}{| \nabla \phi(t,\bfx )|}, \quad \textrm{for $\bfx\in\Gamma(t)$} 
\end{equation}

The evolution of the interface $\Gamma(t)$ is governed by the following advection equation for the level set function: find 
$\phi: I \times \Omega \rightarrow {\bbR}$ such that
\begin{equation}\label{eq:levelsetadv}
\phi_t+\bfbeta \cdot \nabla \phi=0, \quad \phi(0,\bfx) = \phi_0
\end{equation}

Following~\cite{HLZ16ST}, we find an approximation $\phi_h \in V_{0,h/2}$ of the level set function in the space of piecewise linear continuous functions defined on the mesh $\mcK_{0,h/2}$ obtained by refining $\mcK_{0,h}$ uniformly once. 
We consider a continuous piecewise linear approximation $\Gamma_h$ of $\Gamma$ such that $\Gamma_h \cap K$ is a linear segment for $d=2$ and is a subset of a hyperplane in $\bbR^3$, for each $K \in \mcK_{0,h/2}$. A piecewise constant approximation $\bfn_h$ to the exterior unit normal is computed as the spatial gradient of $\phi_h$ (see equation \eqref{eq:normal}). We assume the following hold at every time $t \in I$:
\begin{equation}\label{eq:geomassumptiona}
\| \phi -\phi_h \|_{L^\infty(\Gamma_h)} \lesssim h^2
\end{equation}
and
\begin{equation}\label{eq:geomassumptionb}
\|\bfn^e - \bfn_h \|_{L^\infty(\Gamma_h)} \lesssim h
\end{equation}
Here $\lesssim$ denotes less or equal up to a positive constant,
$\bfn^e$ is the extension of the exact normal to $\Gamma_h$ by the closest point mapping. These assumptions are consistent with the piecewise linear nature of the discrete interface. The subdomain $\Omega_{h,1}$ is defined as the domain enclosed by $\Gamma_h\cup \partial \Omega$ and $\Omega_{h,2}$ as the domain enclosed by $\Gamma_h$.

As in~\cite{HLZ16ST} we use the Crank--Nicolson scheme in time and piecewise linear continuous finite elements with streamline diffusion stabilization in space
to solve the advection equation~\eqref{eq:levelsetadv}. We obtain the method: find $\phi_h^n \in V_{0,h/2}$ 
such that, for $n = 1,2,\dots,N$, 
\begin{align}\label{eq:advdisc}
&(\frac{\phi_h^{n}}{k_n}+\frac{1}{2}\bfbeta^n \cdot \nabla \phi_h^{n},v^n)_\Omega + (\frac{\phi_h^{n}}{k_n}+ \frac{1}{2} \bfbeta^n \cdot \nabla \phi_h^{n},\tau_{SD} \bfbeta^n\cdot \nabla v^n)_\Omega =
\nonumber \\ &\quad =(\frac{\phi_h^{n-1}}{k_n}-\frac{1}{2}\bfbeta^{n-1}\cdot \nabla\phi_h^{n-1}, v^n+\tau_{SD} \bfbeta^n\cdot \nabla v^n)_\Omega  \quad \forall v^n_h \in V_{0,h/2} 
\end{align}
where the streamline diffusion parameter $\tau_{SD}= 2 ( k_n^{-2} + | \bfbeta |^2 h^{-2})^{-1/2}$. 
To keep the level set function a signed distance function, the reinitialization equation, equation (15) in~\cite{SuFa99}, 
can be solved in the same way as the advection equation in~\eqref{eq:advdisc}.

\textbf{An explicit representation}
We use a set of marker points distributed at equal arclength intervals on the interface and a periodic cubic spline as parametrization of the interface. Thus, given a set of markers $\{\bfx_l \}_{l=1}^M \in \bbR^2$ on the interface we have a parametrization $\bfx(\alpha): [\alpha_1, \alpha_M]\rightarrow \bbR^2$ such that 
\begin{align}
\bfx(\alpha_l)&=\bfx_l, \quad l=1, \cdots, M \\
\bfx^m(\alpha_1)&=\bfx^m(\alpha_M), \quad m=0,1,2
\end{align}
where $\bfx(\alpha)=(X_1(\alpha),  X_2(\alpha))$, $X_i$ is a polynomial of degree less or equal to three in each interval $[\alpha_l, \alpha_{l+1}]$, $l=1, \cdots, M$ and has $C^2$ continuity at $\alpha_1,\cdots,\alpha_{M}$ associated with the marker points $x_1,\cdots, x_{M}$. 

The normal is computed from the parametrization $\bfx(\alpha)=(X_1(\alpha),X_2(\alpha))$ as:
\begin{equation}
n_h(\bfx)=\frac{(-X_2'(\alpha),X_1'(\alpha))}{\sqrt{(X_1'(\alpha))^2+(X_2'(\alpha))^2 }}
\end{equation}
Since the interface is smooth we expect the error (measured in max-norm) in the approximation of the geometry and in the approximation of the normal to be:  
\begin{equation}\label{eq:geomassumptionaSp}
\| \Gamma -\Gamma_h \|_{L^\infty(\Gamma_h)} \lesssim h_{\alpha}^4
\end{equation}
and
\begin{equation}\label{eq:geomassumptionbSp}
\|\bfn^e - \bfn_h \|_{L^\infty(\Gamma_h)} \lesssim h_{\alpha}^3
\end{equation}
Here $h_{\alpha}$ is the distance between the marker points and we choose $h_{\alpha}$ to be proportional to the mesh size $h$ in the background mesh.

To evolve the interface the following ordinary differential equation is solved:
\begin{equation}\label{eq:splineadv}
\frac{d \bfx_l}{dt}=\bfbeta(t,  \bfx_l), \quad \bfx_l(0) = \bfx_l^0,  \quad  l=1, \cdots, M
\end{equation}
At each time step a new spline is interpolated through the advected marker points. To avoid clustering or depletion of marker points either a reinitialization step redistributing the points is needed or one can preserve the equal arclength spacing of the marker points by modifying the tangential velocity in equation~\eqref{eq:splineadv}, see~\cite{HoLoSh94}.

\subsubsection{Assembly of the bilinear forms using quadrature in time}\label{sec:assembly}
To compute the space-time integrals in the variational formulation our strategy is to first use a quadrature rule in time and then for each quadrature point compute the integrals in space. 

Using a quadrature formula in time in the interval $I_n$ with quadrature weights $\omega_m^n$ and quadrature points $t_m^n$, $m=1,\dots,n_m$, where $n_m$ is the 
number of quadrature points, recalling equation \eqref{eq:ansatzSurf} and assuming we use linear elements in time ($q=1$), we can approximate the first term in the bilinear form $A_h^n(u,v)$ by
\begin{equation}
\int_{I_n} (\partial_t u,v)_{\Gammah(t)}\, dt\approx
\sum_{m=1}^{n_m} \omega_m^n\frac{1}{k_n}(u_{s,1},v_{s,0})_{\Gammah(t_m^n)} +
\sum_{m=1}^{n_m} \omega_m^n \frac{t_m^n-t_{n-1}}{k_n^2}(u_{s,1},v_{s,1})_{\Gammah(t_m^n)}
\end{equation}
The other space-time integrals are treated in the same way. 

Note that for large time step sizes it may happen that the interface from one quadrature point in time to the next passes over several elements. Thus, it may happen that there are elements in the active mesh $\mcN_{S,h}^n$ (recall equation \eqref{eq:mcnmcs}) which are not intersected by the interface at any quadrature point $t_m^n$,   $m=1,\dots,n_m$ in time.  However, on the faces of those elements the face stabilization is active and therefore the resulting linear system will not be singular.

Consider the closed Newton-Cotes formulas in Table~\ref{tab:quadrature rule}.
\begin{table}[h]
\begin{tabular}{|l|l|l|l|}
\toprule
 $n_m$ & quadrature points  $t_m^n$ & quadrature  weights $\omega_m^n$  & degree of precision \\  \midrule
 $2$ & $t_1^n=t_{n-1}$, $t_2^n=t_{n}$ &$\omega_1^n=\omega_2^n=\frac{k_n}{2}$  & 1\\
 &  & &\\
$3$ & $t_1^n=t_{n-1}$,  $t_3^n=t_{n}$, $t_2^n=\frac{t_{n-1}+t_n}{2}$ & $\omega_1^n=\omega_3^n=\frac{k_n}{6}$, $\omega_2^n=\frac{4k_n}{6}$ & 3\\
& & & \\
$5$ & $t_1^n=t_{n-1}$, $t_5^n=t_{n}$, $t_2^n=\frac{3t_{n-1}+t_n}{4}$, 
&$\omega_1^n=\omega_5^n=\frac{7k_n}{90}$, $\omega_2^n=\omega_4^n=\frac{32k_n}{90}$,  & 5\\ 
& $t_3^n=\frac{t_{n-1}+t_n}{2}$, $t_4^n=\frac{t_{n-1}+3t_n}{4}$  & $\omega_3^n=\frac{12k_n}{90}$  &  \\
\bottomrule
\end{tabular}
\caption{Closed Newton-Cotes formulas. \label{tab:quadrature rule}}
\end{table}
Each quadrature formula integrates exactly polynomials of degree less than or equal to the quadrature formulas degree of precision.
In~\cite{HLZ16ST} we used linear elements in space and time and studied the first two rules in the table above known as the trapezoidal rule and the Simpson's rule, respectively.  In the numerical examples in the next section we use Simpsons' rule when $q=1$ (linear elements are used in time) and the five point Newton-Cotes formula when $q=2$.  Note that these rules include the endpoints of the time interval $I_n$ and some computations can be reused when passing from one space-time slab to another. Note that the quadrature formulas employ equally spaced points and since we need to compute the discrete surface $\Gamma_h$ at the quadrature points we choose the time step size with which we evolve the interface to be $k_n/(n_m-1)$.

\subsection{Numerical examples}\label{sec:numexpS}
We consider an example similar to the last example in~\cite{DeElRa14}. 
The interface $\Gamma$ is an oscillating ellipse defined by  
the zero level set of the level set function,
\begin{equation}
\phi(x,t)=\frac{x_1^2}{(1+0.25\sin(2 \pi t))^2}+x_2^2-1
\end{equation}
where $x=(x_1,x_2) \in \bbR^2$ or by the parametric equations
\begin{align}
X_1(t,\alpha)&=a(t) \cos(\alpha) \nonumber \\ 
X_2(t,\alpha)&=\sin(\alpha)
\end{align}
where $a(t)={1+0.25\sin(2 \pi t)}$. The velocity field is given by
\begin{equation}
\beta=\frac{\pi}{2}\frac{\cos(2\pi t)}{a(t)} (x_1,0)
\end{equation}
The interfacial diffusion coefficient is set to one, $k_S=1$.  We study two different solutions to the surface PDE given in equation~\eqref{eq:uS}, see Example 1 and 2. For p=1 results using the level set method coincide with results using the explicit representation of the interface and we therefore only show results when the interface is represented by a set of markers and a cubic spline parametrization, see Section ~\ref{sec:Intrep}. We always take a large enough number of marker points so that the geometrical error is not dominating the total error.

The computational domain is $[-1.5, \ 1.5]\times [-1.5, \ 1.5]$.  We use a uniform underlying mesh $\mcK_{0,h}$ consisting of triangles with $h=h_{x_1}=h_{x_2}$ and a time step size $k = h/12$.  The error is measured at time $t=0.25$ both in the L$^2$-norm, 
\begin{equation}
\| u^e - u_{h} \|_{L^2(\Gammah(0.25))}
\end{equation}  
and  the $H^1$-norm, 
\begin{equation} 
\| u^e - u_{h} \|_{H^1(\Gammah(0.25))}=\left(\| u^e - u_{h} \|_{L^2(\Gammah(0.25))}^2 + \| (\nablas u)^e -\nablash u_{h} \|_{L^2(\Gammah(0.25))}^2 \right)^{1/2}
\end{equation}
where $u^e$ is an extension of $u$ to $\Gammah$.

\subsubsection{Example 1}
\begin{figure}\centering
\includegraphics[scale=0.45]{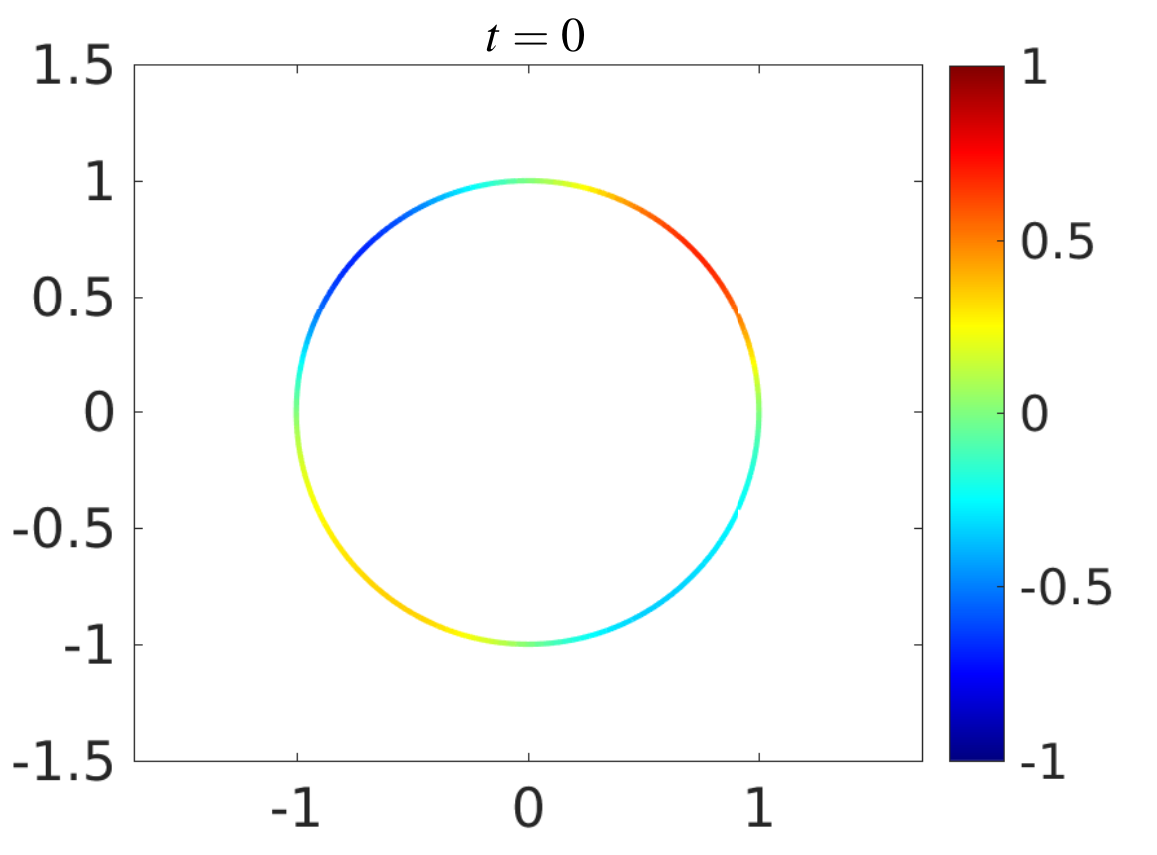}
\includegraphics[scale=0.45]{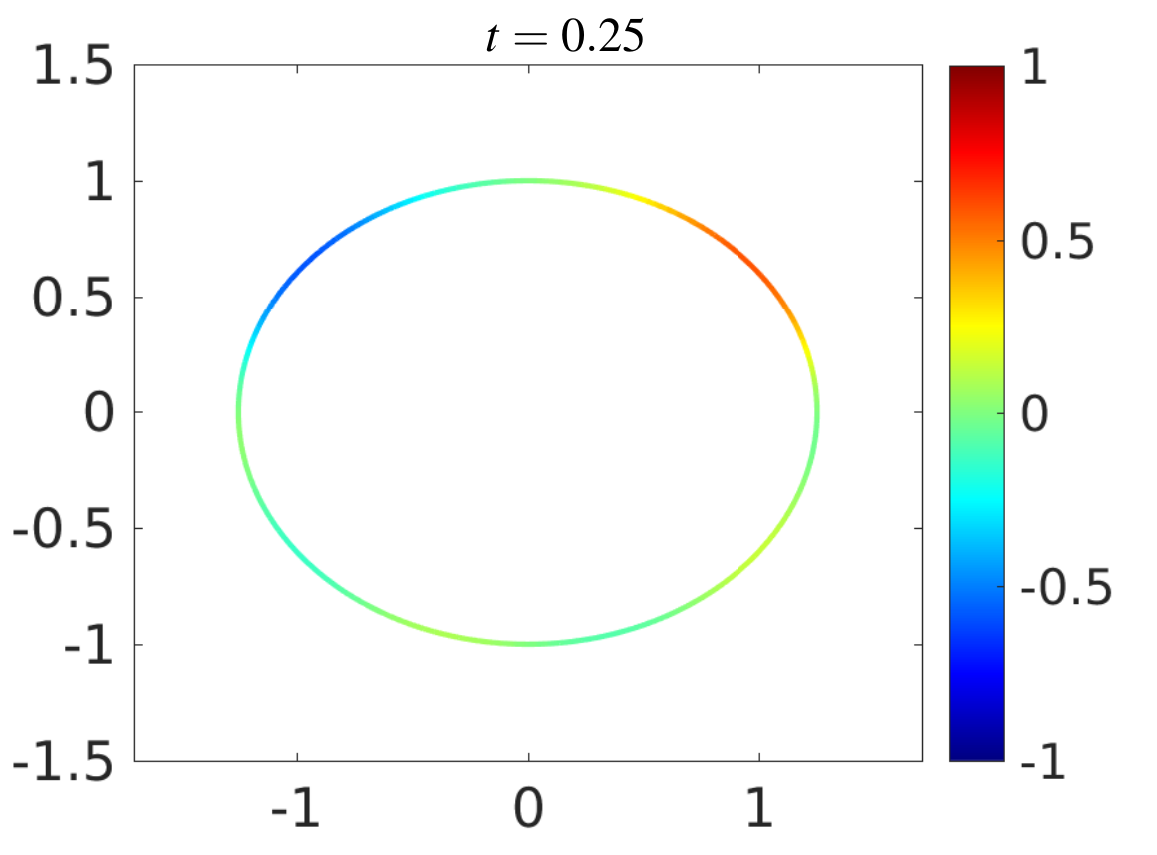} 
\caption{Example 1: The solution of the surface PDE given in equation~\eqref{eq:uS} at times $t=0$ and $t=0.25$ using the proposed space-time CutFEM. The mesh size is $h=0.075$ and time step size is $k = h/12$. \label{fig:solEx1}}.
\end{figure}

As an exact solution of equation~\eqref{eq:uS} we take 
\begin{equation} \label{eq:sol2}
u(t,x)=e^{-4t}x_1x_2+x_1^3x_2^2
\end{equation} 
A right-hand side $f$ to equation~\eqref{eq:uS} is calculated so that the given function~\eqref{eq:sol2} satisfies the surface PDE.  In Fig.~\ref{fig:solEx1} we show the computed solution of equation~\eqref{eq:uS} using the proposed space-time CutFEM with $p=q=2$. The mesh size is $h=0.075$. 
\begin{figure}\centering
\includegraphics[scale=0.6]{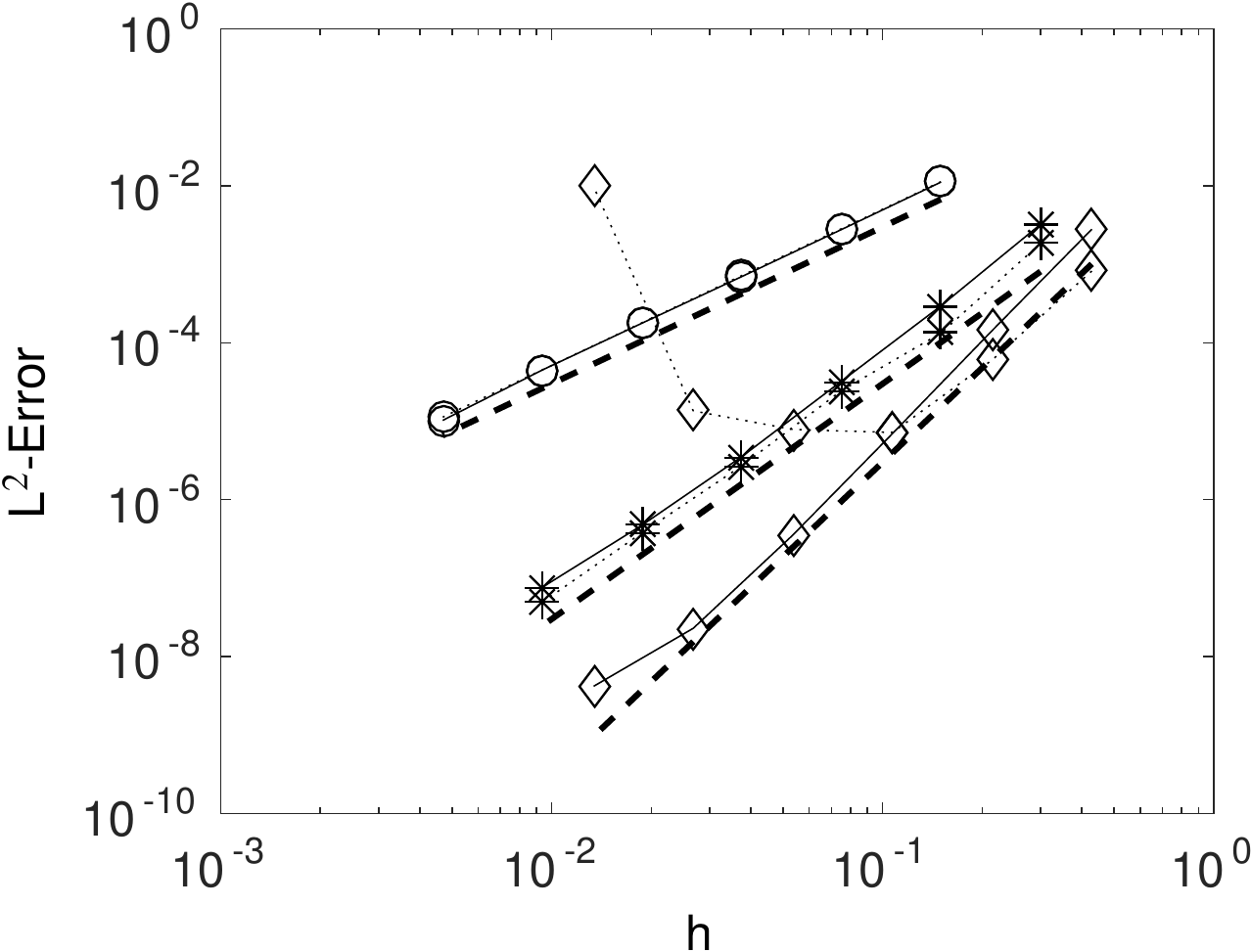} \hspace{2cm}
\includegraphics[scale=0.6]{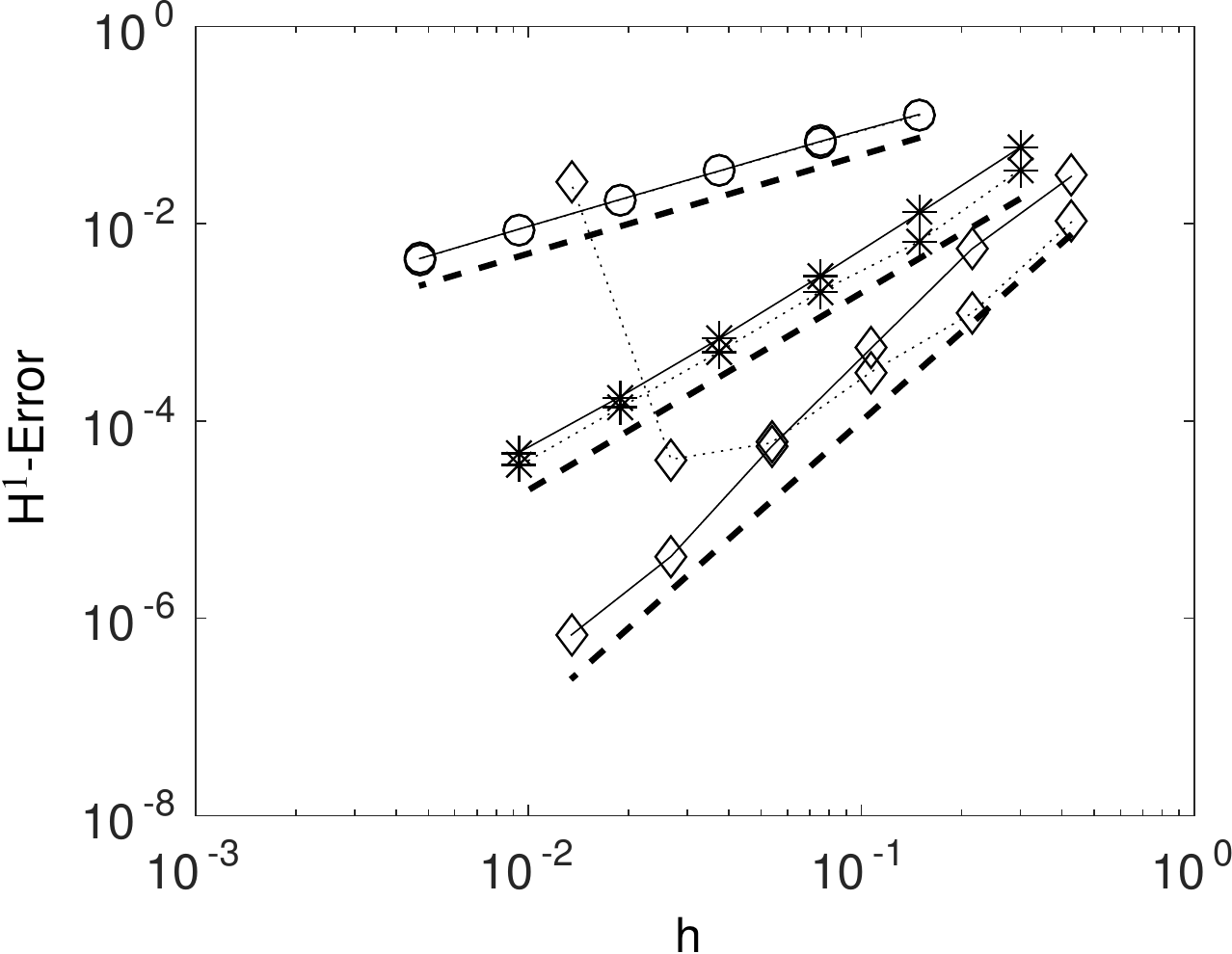} 
\includegraphics[scale=0.6]{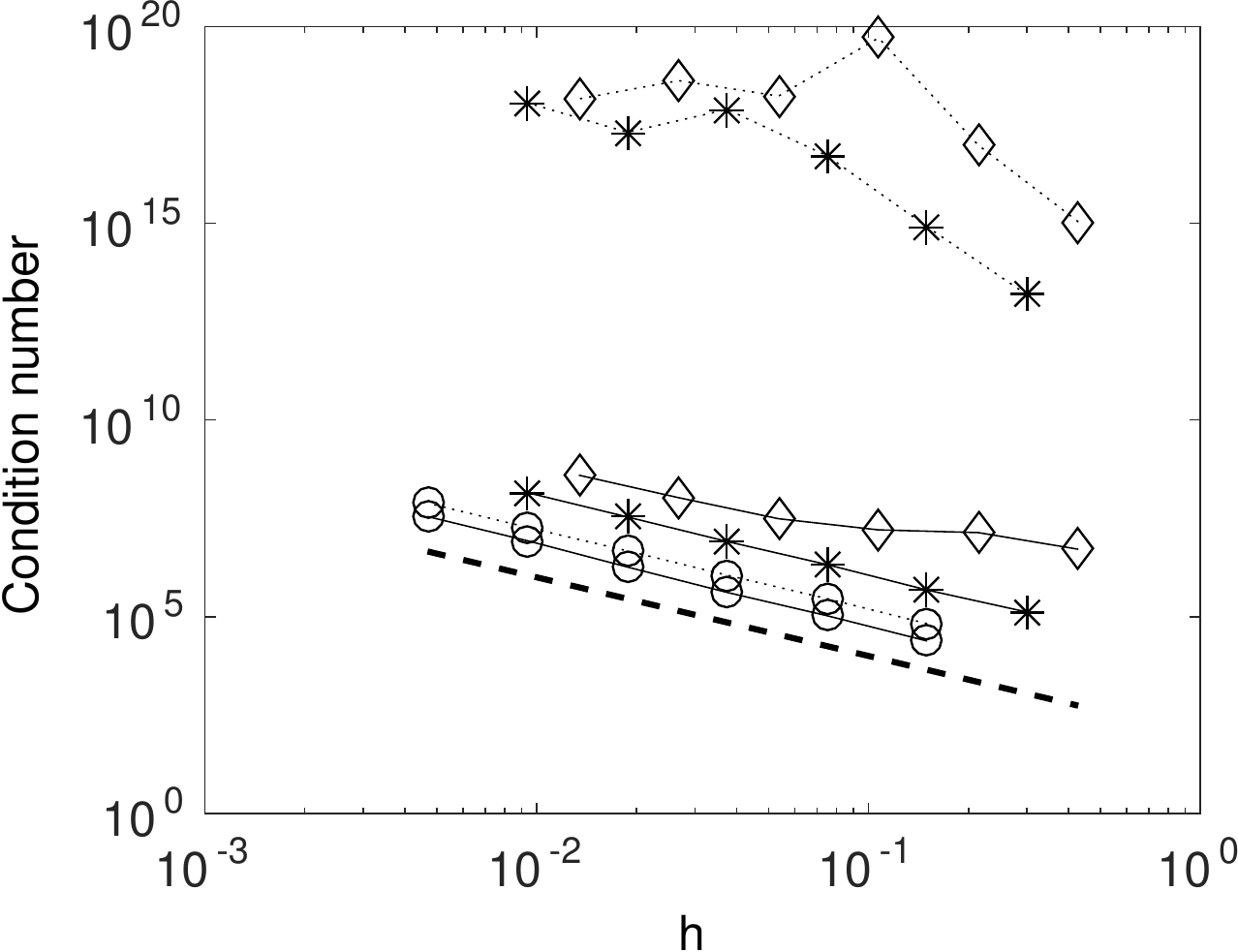} 
\caption{Example 1: The error and condition number versus mesh size $h$ for different degrees of polynomials in the space discretization. Circles: $p=1$. Stars: $p=2$, Diamonds: $p=3$. In time linear elements are used, i.e. $q=1$.
The time step size is $k = h/12$. Symbols connected with a solid line represent results with the proposed stabilization. Symbols connected with a dotted line represent results with only the face stabilization and $\gamma=2i-2$. Top: The error measured in the L$^2$-norm versus mesh size $h$. The dashed lines are proportional to $h^{p+1}$.  Middle: The error measured in the H$^1$-norm versus mesh size $h$.  The dashed lines indicating the expected rates of convergence are proportional to $h^{p}$.  Bottom: The spectral condition number versus mesh size $h$. The dashed line is proportional to $h^{-2}$.  \label{fig:errorvshEx1}}
\end{figure} 

We compare the error in the cut finite element approximation and the condition number of the algebraic system using the new stabilization term, choosing $c_{F,i}=c_{\Gamma,i}=\frac{10^{-1}}{i!}$, with using only the face stabilization term with $\gamma=2i-2$ and $c_{F,i}=\frac{10^{-2}}{i!}$. 
In this example the error in the space discretization dominates and therefore we only show results using linear elements in time.  In space we use linear, quadratic, and cubic elements.  In Fig.~\ref{fig:errorvshEx1} we see the error in the computed solution and the spectral condition number as a function of mesh size $h$. For linear ($p=1$) and quadratic ($p=2$) elements in space the error using the new stabilization and the pure face stabilization with $\gamma=2i-2$ are very simliar and the solid lines and the dotted lines in the figures showing the error in the L$^2$-norm and the H$^1$-norm almost coincide. However, when $p=2$ or $3$ the condition number is very large if only the face stabilization is used. For the cubic elements $p=3$, due to the high condition number the error is dominated by roundoff errors and the convergence in the L$^2$-norm stops and errors increase as the mesh is refined. Diagonal scaling did not improve the condition number. These results show that the face stabilization we used in the space-time CutFEM in~\cite{HLZ16ST} does not control the condition number using higher order elements but the new stabilization term does and the condition number scales as $\mathcal{O}(h^{-2})$ as for standard finite element methods.

Next we consider an example where for $p\geq2$ the error is dominated by the error in the time discretization. 

\subsubsection{Example 2}
\begin{figure}\centering
\includegraphics[scale=0.45]{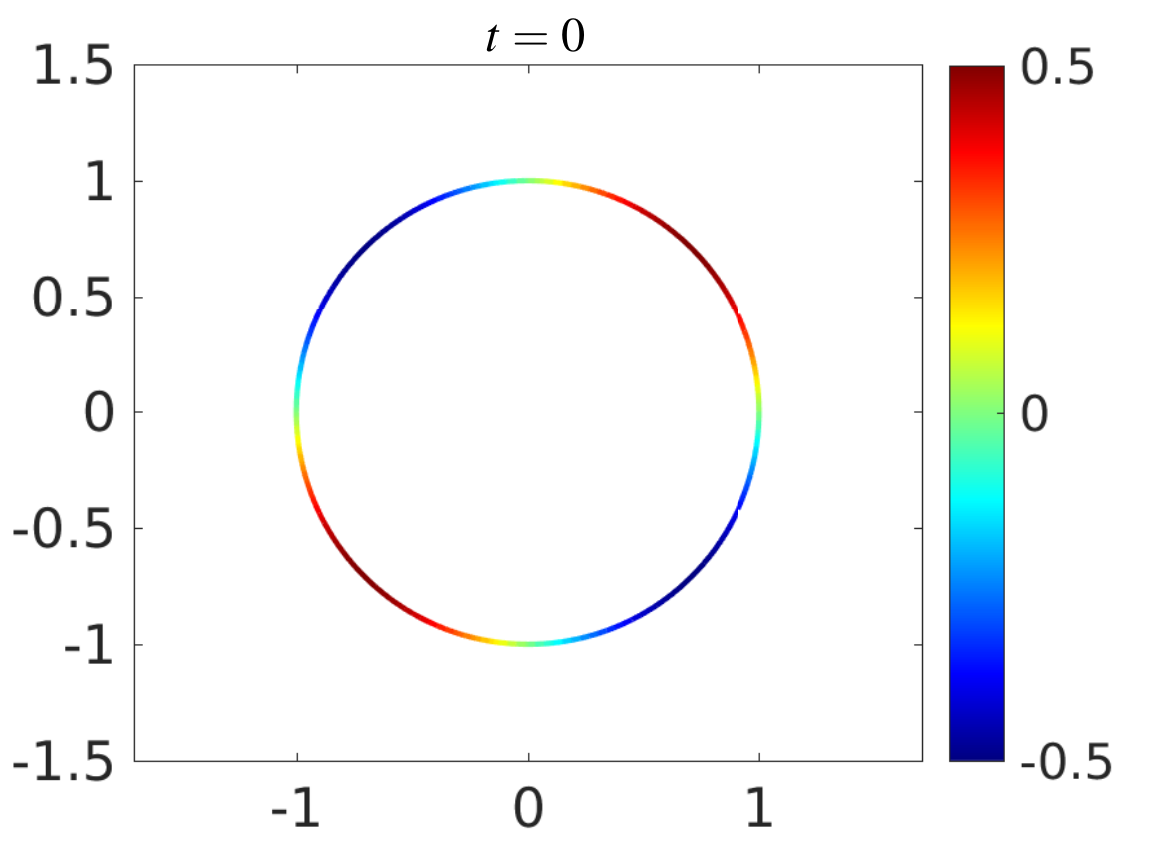} 
\includegraphics[scale=0.45]{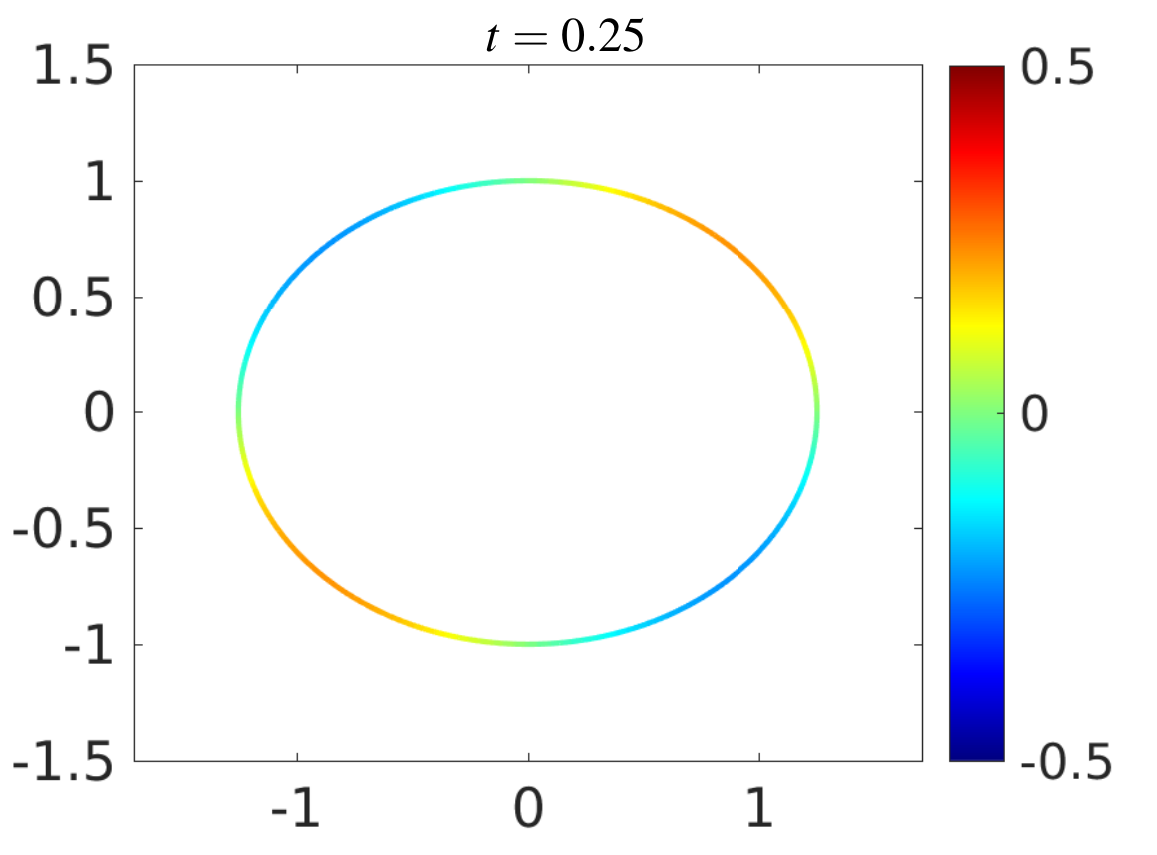} 
\caption{Example 2: The solution of the surface PDE given in equation~\eqref{eq:uS} at times $t=0$ and $t=0.25$ using the proposed space-time CutFEM with mesh size $h=0.075$ and $k = h/12$. \label{fig:solEx2}}
\end{figure}

As an exact solution of equation~\eqref{eq:uS} we now take 
\begin{equation} \label{eq:sol1}
u(t,x)=e^{-4t}x_1x_2
\end{equation}
A right-hand side $f$ to equation~\eqref{eq:uS} is calculated so that $u(x,t)=e^{-4t}x_1x_2$ satisfies the equation. In Fig.~\ref{fig:solEx2} we show the computed solution of equation~\eqref{eq:uS} using the proposed space-time CutFEM with $p=q=2$. Note that p and q are the degree of the polynomials used in space and time, respectively. The mesh size is $h=0.075$.
\begin{figure}\centering
\includegraphics[scale=0.6]{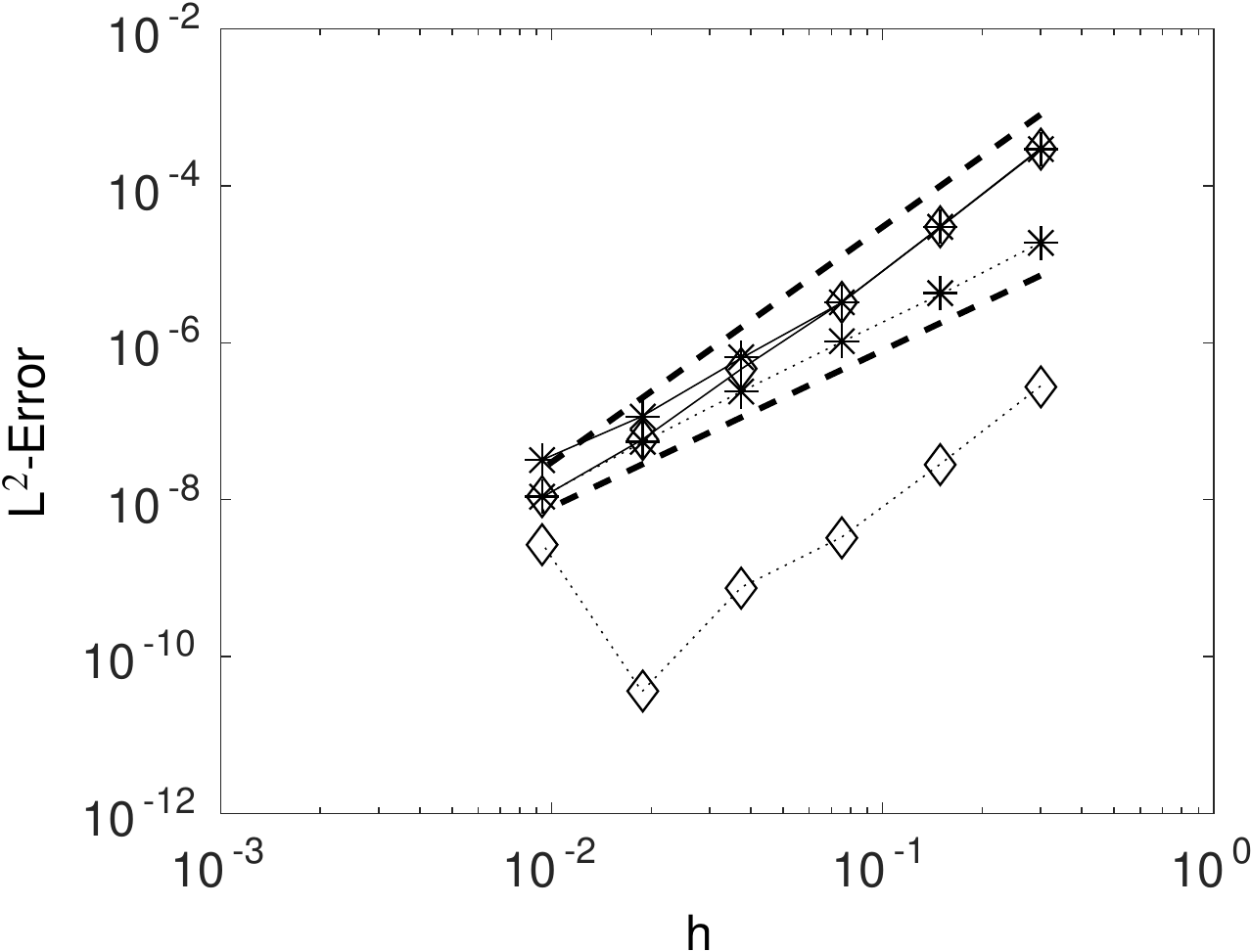} \hspace{2cm}
\includegraphics[scale=0.6]{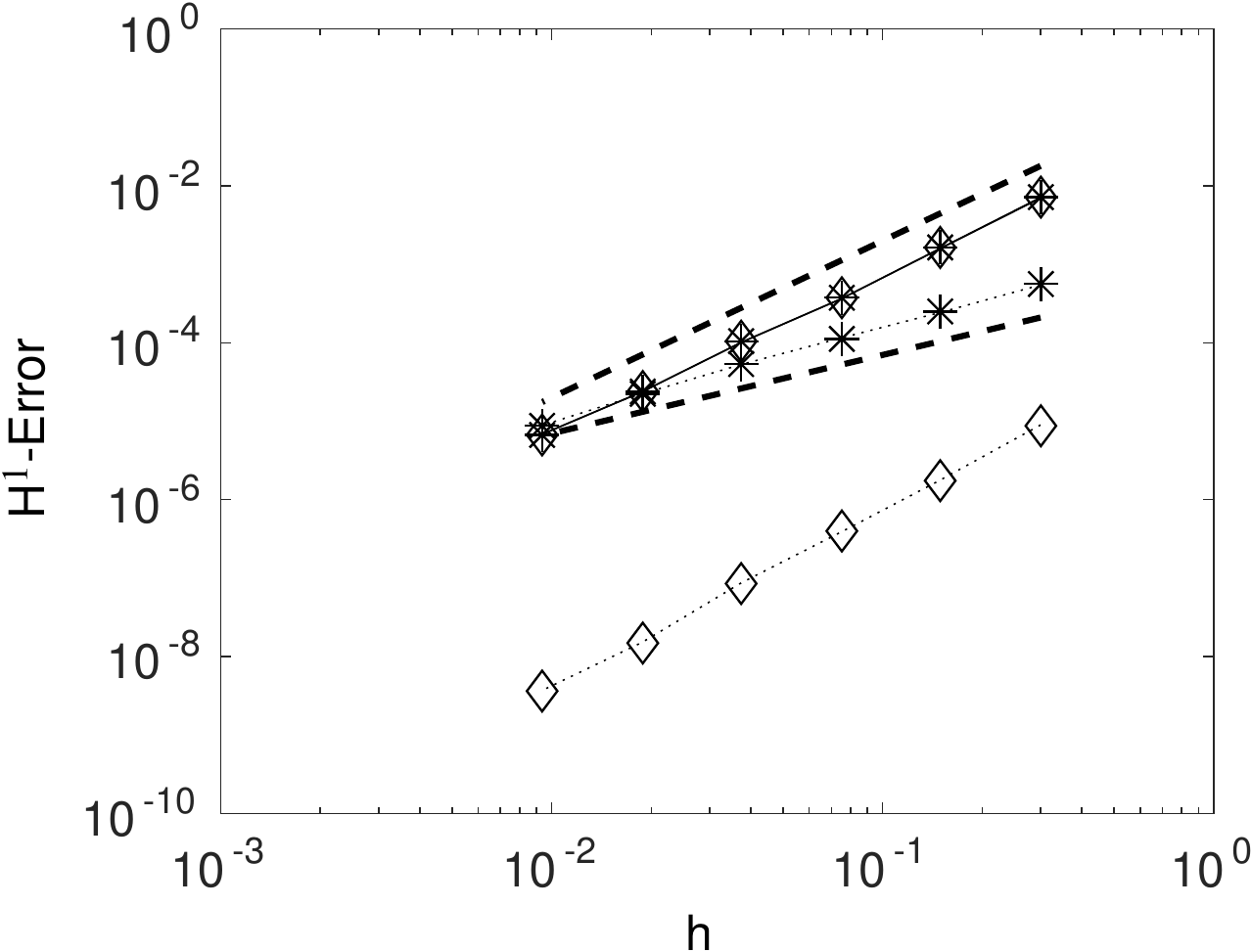} 
\includegraphics[scale=0.6]{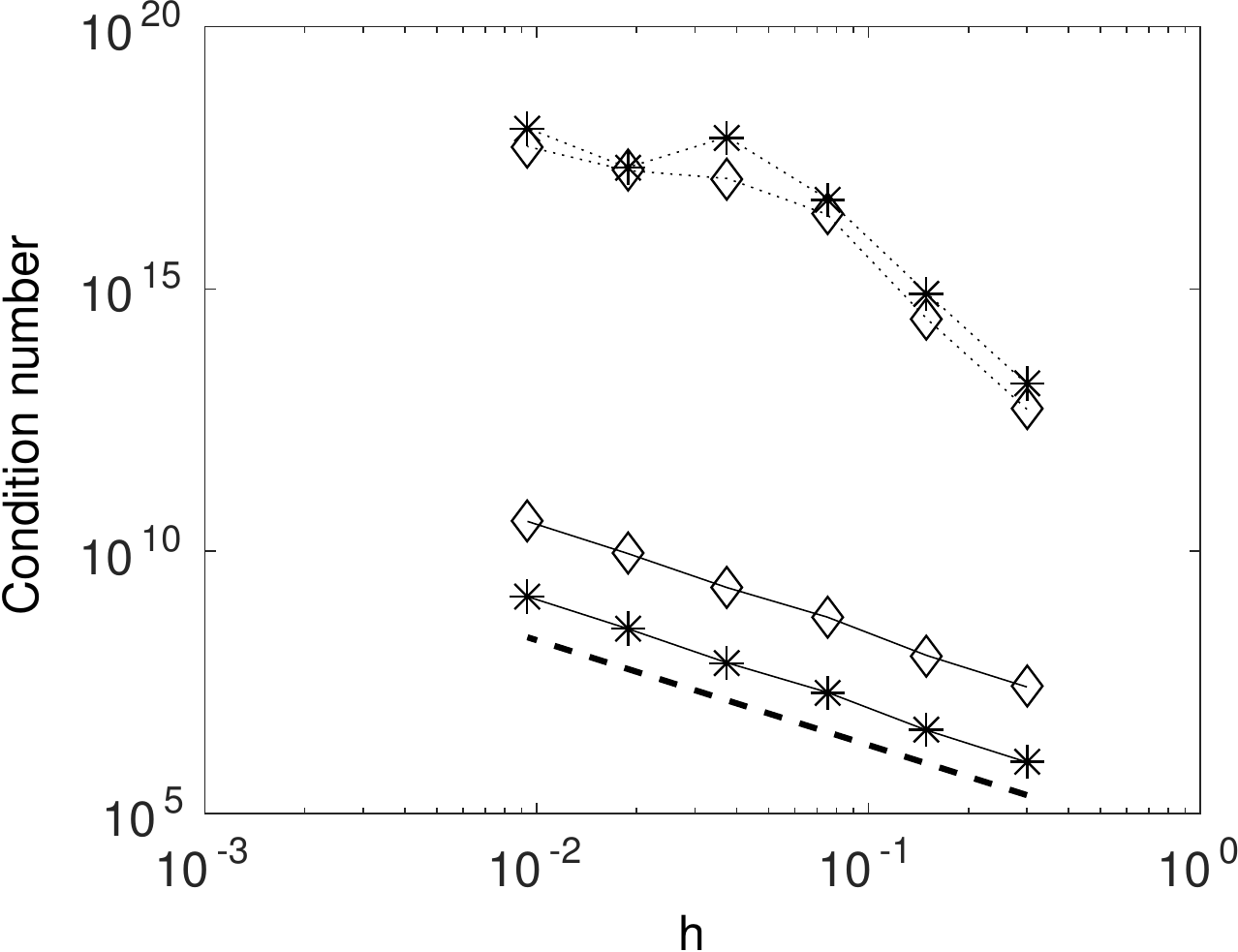} 
\caption{Example 2: The error and condition number versus mesh size $h$ for different degrees of polynomials in the time discretization. Stars: $q=1$, Diamonds: $q=2$. Quadratic elements are used in space, i.e. $p=2$.  The time step size is $k = h/12$. Symbols connected with a solid line represent results with the proposed stabilization. Symbols connected with a dotted line represent results with only the face stabilization and $\gamma=2i-2$.  
Top: The error measured in the L$^2$-norm versus mesh size $h$. The dashed lines are proportional to $h^{q+1}$.  Middle: The error measured in the H$^1$-norm versus mesh size $h$.  The dashed lines are proportional to $h^q$.  Bottom: The spectral condition number versus mesh size $h$. The dashed line is proportional to $h^{-2}$. \label{fig:errorvshEx2}}
\end{figure}

In this example the term $x_1x_2$ in the exact solution is in our finite element space when $p\geq 2$. We study the error and the condition number for $p=2$ with $q=1$ and $q=2$ using the new stabilization term with $c_{F,i}=c_{\Gamma,i}=\frac{10^{-2}}{i!}$ and compare the results with the results using only the face stabilization term with $\gamma=2i-2$ and $c_{F,i}=\frac{10^{-2}}{i!}$. In Fig.~\ref{fig:errorvshEx2} we see that when only the face stabilization is used as in~\cite{HLZ16ST} the error in the time discretization dominates and therefore a big improvement is obtained by using higher order elements in time. However, the condition number is large and the error in the $L^2$-norm starts to increase and is dominated by roundoff errors. A diagonal preconditioning did not improve the condition number. The new stabilization, on the other hand, gives control of the condition number and the spectral condition number increases as $\mathcal{O}(h^{-2})$.  A diagonal preconditioning can now be used to further decrease the condition number. However, with the new stabilization the error in the time discretization is not the dominating error anymore and the errors are not reduced by taking higher order elements in time. We obtain similar results using cubic elements in space.

In Fig.~\ref{fig:errorvshEx2} we see second order convergence of the error in the $L^2$-norm when $p=2$ and $q=1$, regardless of which of the two stabilization terms we used. For coarse meshes the error of the cut finite element approximation obtained using the new stabilization term converges faster. We obtain third order convergence when we use quadratic elements in both space and time, i.e., $p=q=2$. Using the face stabilization we obtain third order convergence initially but the convergence stops when the condition number becomes too large.

\subsubsection{Discussion}
In the L$^2$-norm convergence orders $p+1$ in space and $q+1$ in time has been observed. For discontinuous Galerkin methods in time based on polynomials of order $q$ superconvergence, i.e. convergence of order $2q+1$ in the nodes $t_n$ has been reported, see \cite{Th06}. We did not observe such superconvergence. With the new stabilization term the condition number scales as $\mathcal{O}(h^{-2})$ independent of the order of the elements we use, as in standard finite element methods. Using only the face stabilization however resulted in large condition numbers that sometimes increased as $\mathcal{O}(h^{-6})$. In the last example the new stabilization term resulted in large errors compared to using only the face stabilization. However, we emphasize that the second example is a very special case since the term $x_1x_2$ in the exact solution is in our finite element space when p=2 and we therefore see the error from the stabilization term. 
In summary, we observe that the new stabilization is strong enough to control the condition number also for higher order elements and weak enough to not destroy the convergence order of the method.

\section{A coupled bulk-surface problem}\label{sec:modelcoupled}
We now consider a coupled bulk-surface problem modeling the evolution of soluble surfactants. In non-dimensional form we have
\begin{alignat}{2}
\partial_t u_B + \bfbeta \cdot \nabla u_B  - \nabla \cdot \left(\frac{1}{Pe} \nabla u_B \right) &= 0  
\quad &&\text{in $\Omega_1(t)$}   \label{eq:uBPN} \\
-\bfn \cdot \frac{1}{\textrm{Pe}} \nabla u_B &= \textrm{Da}f_{\textrm{coupling}}  
\quad &&\text{on $\Gamma(t)$} \label{eq:BC1N} \\ 
-\bfn_{\partial \Omega} \cdot \frac{1}{\textrm{Pe}} \nabla u_B &= 0 
\quad &&\text{on $\partial \Omega$} \label{eq:BC2N} \\
\partial_t u_S + \bfbeta \cdot \nabla u_S +(\divs \bfbeta) u_S  -\frac{1}{\textrm{Pe}_{\textrm{S}}} \Delta_\Gamma u_S &=   f_{\textrm{coupling}}
\quad &&\text{on $\Gamma(t)$} \label{eq:uSPN}
\end{alignat}
for all $t\in I$ with 
\begin{equation}\label{eq:coupnondim} 
f_{\textrm{coupling}}=\alpha u_B(1-u_S)-\textrm{Bi}u_S
\end{equation} 
given from the Langmuir model. Examples of other models can be found in for example~\cite{RaFeLi00}. The non-dimensional numbers Pe and Pe$_{\textrm{S}}$ are the bulk and interfacial Peclet numbers, 
Da is the Damk\"{o}hler number, Bi is the Biot number, and
$\alpha=\frac{k_aLu_B^\infty}{\beta^\infty}$, where $k_a$ is the adsorption coefficient, $L$, $\beta^\infty$, and $u_B^\infty$ are the characteristic values for length, velocity, and bulk surfactant concentration~\cite{GaTo12}.  
The conservation of surfactants is expressed in non-dimensional form as:
\begin{equation}\label{eq:conservN}
\int_{\Omega_1(t)} u_B dv +\textrm{Da}\int_{\Gamma(t)} u_S ds=\overline{u}_0 \quad \forall t \in I
\end{equation}
Initial conditions $u_B(0,\bfx) = u_B^0$ in $\Omega_1(0)$ and $u_S(0,\bfx) = u_S^0$ on $\Gamma(0)$ are given. 
Note that the surfactant is soluble only in the outer fluid phase $\Omega_1(t)$, this is not a restriction of the method but a simplification.

\subsection{The space-time cut finite element method}\label{sec:method}
We now use the space-time cut finite element method presented in~\cite{HLZ16ST} 
with linear elements in both space and time (i.e. $p=q=1$). We again use discontinuous elements in time and solve the discrete equations one space-time slab at a time and in the time interval $I_n$ the solution throughout the current slab will depend only on the solution at $t_{n-1}^-$.
We follow~\cite{HLZ16ST}. 

\subsubsection{Mesh and spaces}\label{sec:meshandspacec}
\begin{figure}\centering
\includegraphics[scale=0.45]{IllustactivemeshS} \hspace{0.5cm}
\includegraphics[scale=0.45]{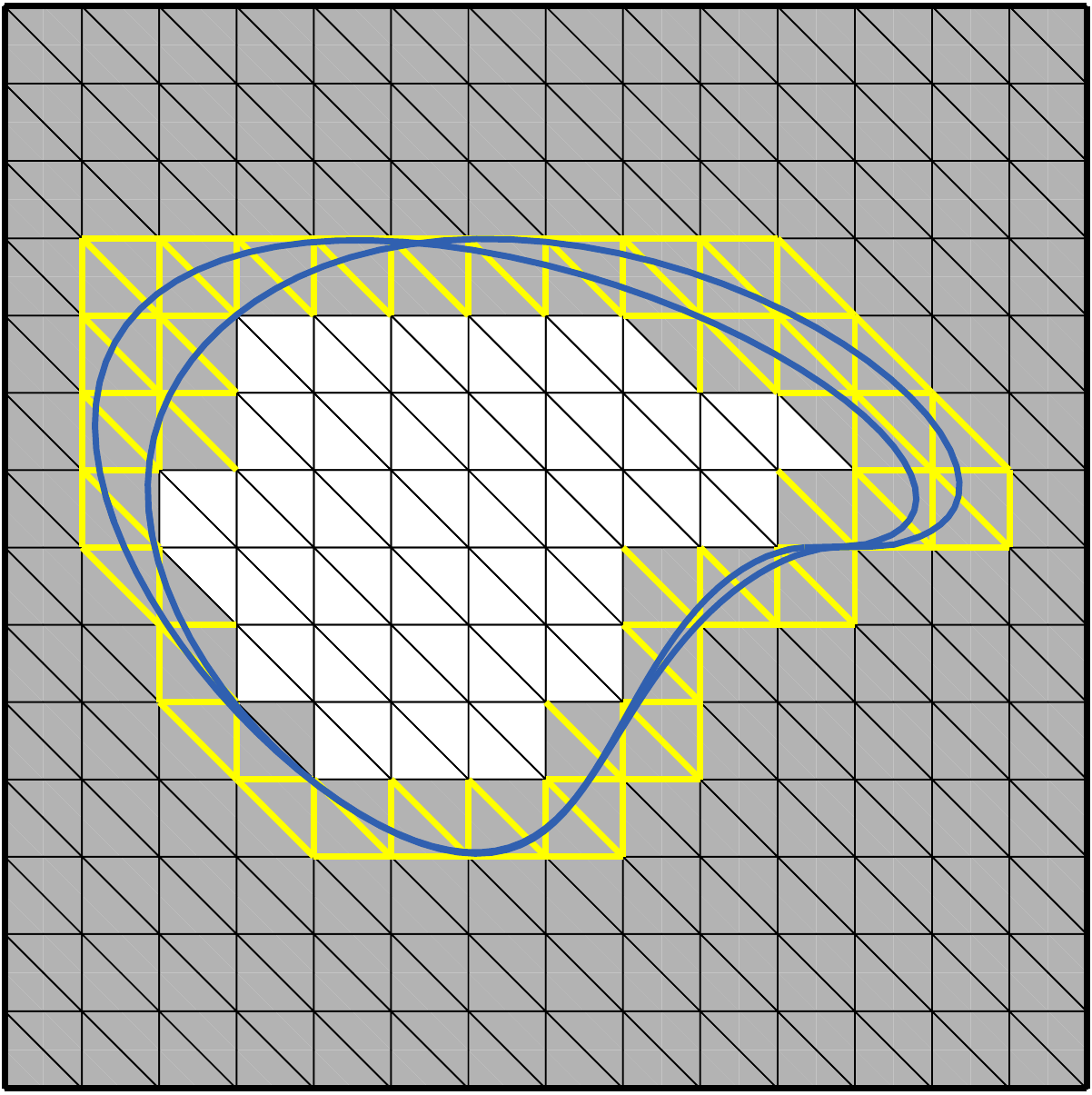}
\caption{Illustration of the sets introduced in Section~\ref{sec:meshandspacec}. In both figures the two blue curves show the position of the interface at the endpoints $t=t_{n-1}$ and $t=t_{n}$ of the time interval $I_n = (t_{n-1},t_n]$.  Left:  the shaded domain shows $\mcN_{S,h}^n$ and edges in $\mcF_{S,h}^n$  are marked with yellow thick lines.  Right: the shaded domain shows  $\mcN_{B,h}^n$ and edges in $\mcF_{B,h}$ are marked with yellow thick lines.  \label{fig:illust}}.
\end{figure}

Define the following sets 
\begin{equation}
\mcK_{B,h}(t) = \{K \in \mcK_{0,h} \,:\, K \cap \Omegaoh(t) \neq \emptyset \},
\,  \
\mcK_{S,h}(t) = \{K \in \mcK_{0,h} \,:\, K \cap \Gammah(t) \neq \emptyset \}
\end{equation}
and the active meshes 
\begin{equation}\label{eq:mcnmcscoup}
\mcN_{B,h}^n = \bigcup_{t \in I_n} \bigcup_{K \in \mcK_{B,h}(t)} K,
\qquad
\mcN_{S,h}^n = \bigcup_{t \in I_n} \bigcup_{K \in \mcK_{S,h}(t)} K  
\end{equation}
As in Section~\ref{sec:meshspaceSurf} $\mcK_{0,h}$ is the fixed background mesh 
and  $I_n = (t_{n-1},t_n]$ is of length $k_n = t_n - t_{n-1}$ for $n = 1,2,\dots,N$. The active meshes are illustrated in Fig.~\ref{fig:illust} by the shaded domains. 

Associated to the active meshes are the space-time slabs
$S_{B}^n = I_n \times \mcN_{B,h}^n$ 
and
$S_S^n = I_n \times \mcN_{S,h}^n$ 
on  which we define the spaces
\begin{equation}
V_{B,h}^n = P_1(I_n) \otimes V_{0,h}|_{\mcN_{B,h}^n},
\qquad V_{S,h}^n=P_1(I_n) \otimes V_{0,h}|_{\mcN_{S,h}^n}
\end{equation}
where $V_{0,h}$ is the space of continuous piecewise linear polynomials defined on the background mesh $\mcK_{0,h}$,
and we let 
\begin{equation}
W_h^n = V_{B,h}^n \times V_{S,h}^n
\end{equation}

Functions in $W_h^n$ take the form
\begin{equation}\label{eq:ansatz}
v_h(t,\bfx) = (v_{B,h},v_{S,h})=\left(v_{B,0}+v_{B,1}\frac{t-t_{n-1}}{k_n},v_{S,0}+v_{S,1}\frac{t-t_{n-1}}{k_n} \right)
\end{equation}
where $t\in I_n$ and $v_{B,j}$ and $v_{S,j}$, $j=0,1$ can be written as 
\begin{align}\label{eq:funcwh}
v_{B,j}&= \sum_{i=1}^{N_B} \xi_{ij}^B \varphi_i(\bfx)|_{\mcN_{B,h}^n}, \quad 
v_{S,j}= \sum_{i=1}^{N_S}  \xi_{ij}^S \varphi_i(\bfx)|_{\mcN_{S,h}^n} 
\end{align}
Here $\xi_{ij}^B, \xi_{ij}^S \in \bbR$ are coefficients, $\varphi_i(\bfx)$ is the standard nodal basis function associated with mesh vertex $i$, $N_B$ and $N_S$ are the number of nodes in $\mcN_{B,h}^n$ and in $\mcN_{S,h}^n$, respectively. 

\subsubsection{The variational formulation} 
Assuming $\alpha$, Bi, Da are positive constants multiplying the bulk PDE, equation \eqref{eq:uBPN} with a test function $\frac{\alpha}{\textrm{Da}} v_B$ and the surface PDE, equation \eqref{eq:uSPN} with a test function $\textrm{Bi} v_s$, integrating by parts, and using the boundary conditions, equation \eqref{eq:BC1N}-\eqref{eq:BC2N}, yields the weak form.  Given $u_h(t_{n-1}^-,\bfx)$ and $\overline{u}_0$ (see equation~\eqref{eq:conservN}) we consider the following weak formulation: find $u_h=(u_B,u_S)\in W_h^n$ and $\lambda \in \bbR$, such that 
\begin{align}\label{eq:spacetimeform}
A_h^n(u_h,v_h) + J_h^n (u_h,v_h) 
&  + \lambda\left( (1,v_B)_{\Omega_{h,1}(t_n)} +\textrm{Da}(1,v_S)_{\Gammah(t_n)} \right) 
\nonumber \\
& + \mu \left((u_B,1)_{\Omega_{h,1}(t_n)}+\textrm{Da} (u_S,1)_{\Gammah(t_n)}\right) 
=\mu \overline{u}_0,  
\end{align}
for all $v_h \in W_h^n, \mu \in \bbR$. Here 
\begin{align}
A_h^n (u,v) &= \int_{I_n} \frac{\alpha}{\textrm{Da}}(\partial_t u_B,v_B)_{\Omega_{h,1}(t)}\, dt+ \int_{I_n} \textrm{Bi}(\partial_t u_S,v_S)_{\Gammah(t)} \, dt + \int_{I_n} a_h(t,u,v)\, dt 
\nonumber \\
&\quad 
-\int_{I_n} \alpha (u_Bu_S, \alpha v_B- \textrm{Bi} v_S)_{\Gammah(t)} \, dt 
\nonumber \\
&\quad 
+\frac{\alpha}{\textrm{Da}}([u_B],v_B(t_{n-1}^+,\bfx))_{\Omega_{h,1}(t_{n-1})} 
+\textrm{Bi}([u_S],v_S(t_{n-1}^+,\bfx))_{\Gammah(t_{n-1})} 
\end{align}
with
\begin{equation}\label{{eq:defa}}
a_h(t,u,v)=\frac{\alpha}{\textrm{Da}} a_{B,h}(t,u_B,v_B)+\textrm{Bi}a_{S,h}(t,u_S,v_S)+a_{BS,h}(t,u,v)
\end{equation}
and
\begin{equation}\label{eq:contformah}
\begin{cases}
a_{B,h}(t,u_B,v_B)= (\bfbeta \cdot \nabla u_B, v_B)_{\Omega_{h,1}(t)}  + \left(\frac{1}{Pe} \nabla u_B,\nabla v_B\right)_{\Omega_{h,1}(t)} 
\\ 
a_{S,h}(t,u_S,v_S)= (\bfbeta \cdot \nabla u_S,v_S)_{\Gammah(t)} +((\divsh \bfbeta) u_S,v_S)_{\Gammah(t)}  + \left(\frac{1}{Pe_{\textrm{s}}} \nablash u_S,\nablash v_S \right)_{\Gammah(t)}
 \\
a_{BS,h}(t,u,v)=(\alpha u_B - \textrm{Bi} u_S, \alpha v_B - \textrm{Bi} v_S)_{\Gammah(t)} 
\end{cases}
\end{equation}

To stabilize the method we use a face stabilization of the form 
\begin{equation}
J_h^n(u_h,v_h)= \int_{I_n}  \tau_B h j_B(u_{B,h},v_{B,h}) + \tau_S  j_S (u_{S,h},v_{S,h})\, dt
\end{equation}
where $\tau_B, \tau_S $ are positive parameters and
\begin{align}
j_B(v_B,w_B) &= \sum_{F \in \mcF_{B,h}} ([\bfn_F\cdot \nabla v_B],[\bfn_F\cdot \nabla w_B])_F
\\
j_S(v_S,w_S) &= \sum_{F \in \mcF_{S,h}} ([\bfn_F \cdot \nabla v_S],[\bfn_F\cdot \nabla w_S])_F
\end{align}
Here $\mcF_{S,h}$ is the set of internal faces in the active surface mesh $\mcN_{S,h}^n$ and $\mcF_{B,h}$ is the set of faces that 
are internal in the active bulk mesh $\mcN_{B,h}^n$ and also belong to an element in $\mcN_{S,h}^n$, see Fig.~\ref{fig:illust}.

\subsection{Implementation}\label{sec:implement}
Since the bulk and the surface surfactant concentrations are coupled 
through a nonlinear term, see \eqref{eq:coupnondim}, the proposed method~\eqref{eq:spacetimeform} leads 
to a nonlinear system of equations in each time step, 
which we solve using Newton's method. To formulate 
Newton's method we define the residual $F$ and the 
Jacobian $DF$ as follows 
\begin{align}\label{eq:F}
F(u,\lambda)&= \int_{I_n} \frac{\alpha}{\textrm{Da}} (\partial_t u_B,v_B)_{\Omega_{h,1}(t)}\, dt+ \int_{I_n} \textrm{Bi}(\partial_t u_S,v_S)_{\Gammah(t)} \, dt + \int_{I_n} a_h(t,u,v)\, dt +
\nonumber \\
&
-\int_{I_n}\alpha (u_Bu_S, \alpha v_B-\textrm{Bi}v_S)_{\Gammah(t)} \, dt  
+ \frac{\alpha}{\textrm{Da}}([u_B],v_B(t_{n-1}^+,\bfx))_{\Omega_{h,1}(t_{n-1})} 
\nonumber \\
& 
+\textrm{Bi}([u_S],v_S(t_{n-1}^+,\bfx))_{\Gammah(t_{n-1})} 
+\int_{I_n}j_h(u,v)\, dt  -\mu \overline{u}_0
\nonumber \\
&
+\lambda\left( (1,v_B)_{\Omega_{h,1}(t_n)} +\textrm{Da} (1,v_S)_{\Gammah(t_n)} \right) + \mu \left((u_B,1)_{\Omega_{h,1}(t_n)}+\textrm{Da} (u_S,1)_{\Gammah(t_n)}\right) 
\end{align}
\begin{align}\label{eq:DF}
&DF(u,\lambda)(w,\hat{\lambda}) = \int_{I_n}  \frac{\alpha}{\textrm{Da}}(\partial_t w_B,v_B)_{\Omega_{h,1}(t)}\, dt+ \int_{I_n} \textrm{Bi}(\partial_t w_S,v_S)_{\Gammah(t)} \, dt + \int_{I_n} a_h(t,w,v)\, dt 
\nonumber \\
&
-\int_{I_n}\alpha (w_Bu_S, \alpha v_B-\textrm{Bi}v_S)_{\Gammah(t)} \, dt -\int_{I_n}\alpha (u_Bw_S, \alpha v_B-\textrm{Bi}v_S)_{\Gammah(t)} \, dt 
\nonumber \\
&
+ \frac{\alpha}{\textrm{Da}}(w_B,v_B(t_{n-1}^+,\bfx))_{\Omega_{h,1}(t_{n-1})} 
+\textrm{Bi}(w_S,v_S(t_{n-1}^+,\bfx))_{\Gammah(t_{n-1})} 
+\int_{I_n}j_h(w,v)\, dt 
\nonumber \\
&
+\hat{\lambda}\left( (1,v_B)_{\Omega_{h,1}(t_n)} +\textrm{Da} (1,v_S)_{\Gammah(t_n)} \right) + \mu \left((w_B,1)_{\Omega_{h,1}(t_n)}+\textrm{Da} (w_S,1)_{\Gammah(t_n)}\right)
\end{align}
With this notation the nonlinear problem resulting from (\ref{eq:spacetimeform}) takes the form: 
find $u_h\in W_h^n$ and $\lambda \in \bbR$ such that 
$F(u_h,\lambda)=0$, and the corresponding Newton's method reads: 
\begin{enumerate}
\item[1.] Choose initial guesses $u_{h,0}$ and $\lambda_0$
\item[2.]  while  $||(w,\hat{\lambda})||>$ tol
\begin{itemize}
\item Solve: $DF(u_{h,0},\lambda_0)(w,\hat{\lambda}) =F(u_{h,0},\lambda_0)$
\item Update $u_{h,0}$: $u_{h,0}=u_{h,0}-w$ and $\lambda_0$: $\lambda_0=\lambda_0-\hat{\lambda}$
\end{itemize}
\end{enumerate}
For $t\in I_n$ we choose the initial guess $u_{h,0}$ to be the solution at $t_{n-1}^- $, i.e. $u_{h,0}(t,\bfx)=u_h(t_{n-1}^-,\bfx)$.  
 
As before, we approximate the space-time integrals using Simpson's rule, see Table \ref{tab:quadrature rule}. At each time interval $I_n$, we compute the discrete surface $\Gamma_h$ at the quadrature points $t_m^n$  as the 
zero level set of the approximate signed distance function $\phi_h(t_m^n,x)$. 
The intersection $\Gamma_h(t_m^n) \cap K$ is planar, since $\phi_h$ is piecewise linear, and we can therefore easily compute the contribution of the surface integrals to the stiffness matrix. The contribution from integration on $\Omega_{h,1}(t_q^n) \cap K$ is divided into contributions on one or several triangles in two dimensions and tetrahedra in three dimensions depending on how the interface cuts element $K$.

Finally, we use a direct solver to solve the linear system of $2(N_B+N_S)+1$ equations: 
\begin{equation}
DF(u_{h,0},\lambda_0)(w,\hat{\lambda}) =F(u_{h,0},\lambda_0)
\end{equation}
for $\hat{\lambda} \in \bbR$ and
\begin{equation}
w=\left(\begin{array}{l}
w_{B,0} \nonumber \\
w_{S,0}\nonumber \\
w_{B,1} \nonumber \\
w_{S,1} \nonumber \\
\end{array} \right)
\end{equation}

\subsection{Numerical example}\label{sec:numexp}
We use one of the examples in~\cite{HLZ16ST}.  The coupled bulk-surface problem is from Section 5.3 of~\cite{ChLai14}. The initial interface is a circle with radius $r_0=0.3$ centered in $(x_0, y_0)=(0.1,0)$ and the velocity field is given by 
\begin{equation}
\bfbeta=\left(-\frac{1}{2}(1+\cos(\pi x))\sin(\pi y), \frac{1}{2}(1+\cos(\pi y))\sin(\pi x) \right)
\end{equation} 
The computational domain is chosen as $\Omega=[-1, 1] \times [-1, 1]$. A uniform fixed background mesh $\mcK_{0,h}$ consisting of triangles of size $h$ is used and a constant time step size of the form $k = h/8$. The non-dimensional numbers are set to $\textrm{Pe}=\textrm{Pe}_\textrm{S}=100$ and $\textrm{Bi}=\alpha=\textrm{Da}=1$. The initial surface and bulk surfactant concentrations are $u_S(0,x,y)=0$ and 
\begin{equation}
u_B(0,x,y)=\left\{
\begin{array}{l}
0.5(1-x^2)^2 \quad \textrm{if $r>1.5r_0$} \nonumber \\
0.5(1-x^2)^2w(r) \quad \textrm{if $r_0\leq r \leq 1.5r_0$} \nonumber \\
0 \quad \textrm{otherwise} 
\end{array}\right.
\end{equation}
with $r=\sqrt{(x-x_0)^2+(y-y_0)^2}$ and
\begin{equation}
w(r)=\frac{1}{2}\left(1-\cos \left(\frac{(r-r_0)\pi}{0.5r_0}\right)\right)
\end{equation}
 In the computations the stabilization constants $\tau_B$ and $\tau_S$ are $10^{-2}$.

\begin{figure}
\begin{center}  
\includegraphics[width=0.45\textwidth]{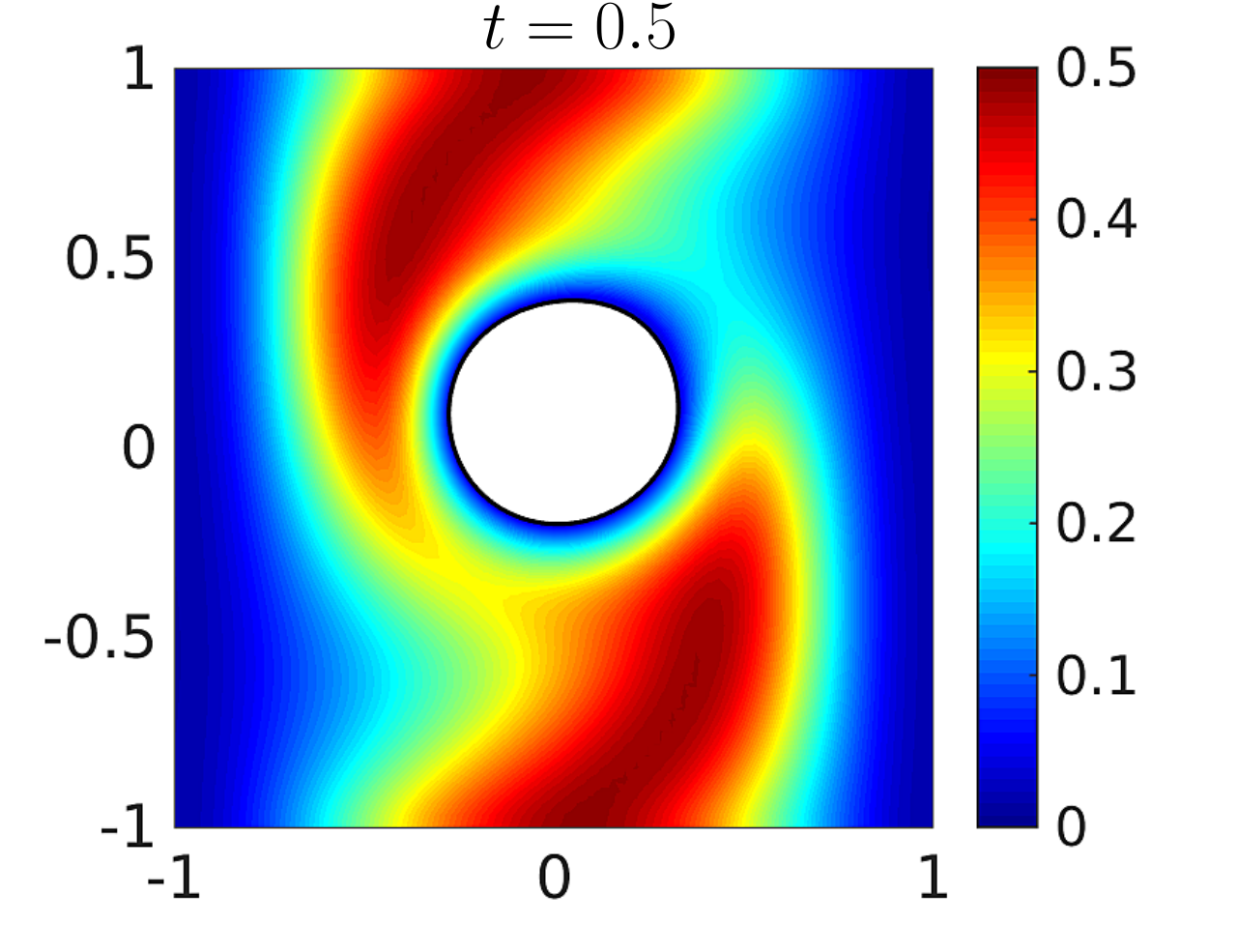}
\includegraphics[width=0.45\textwidth]{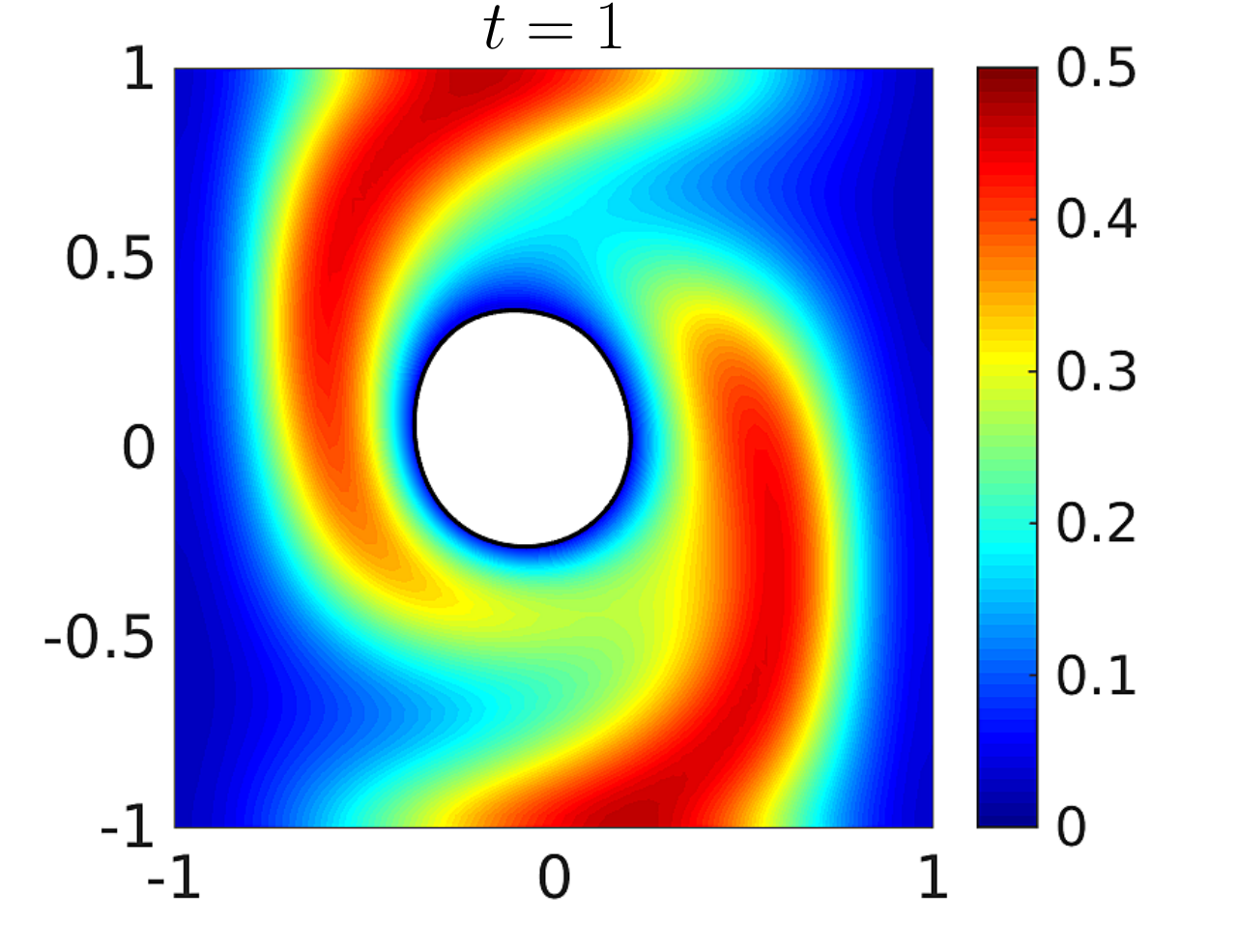}\\
\includegraphics[width=0.45\textwidth]{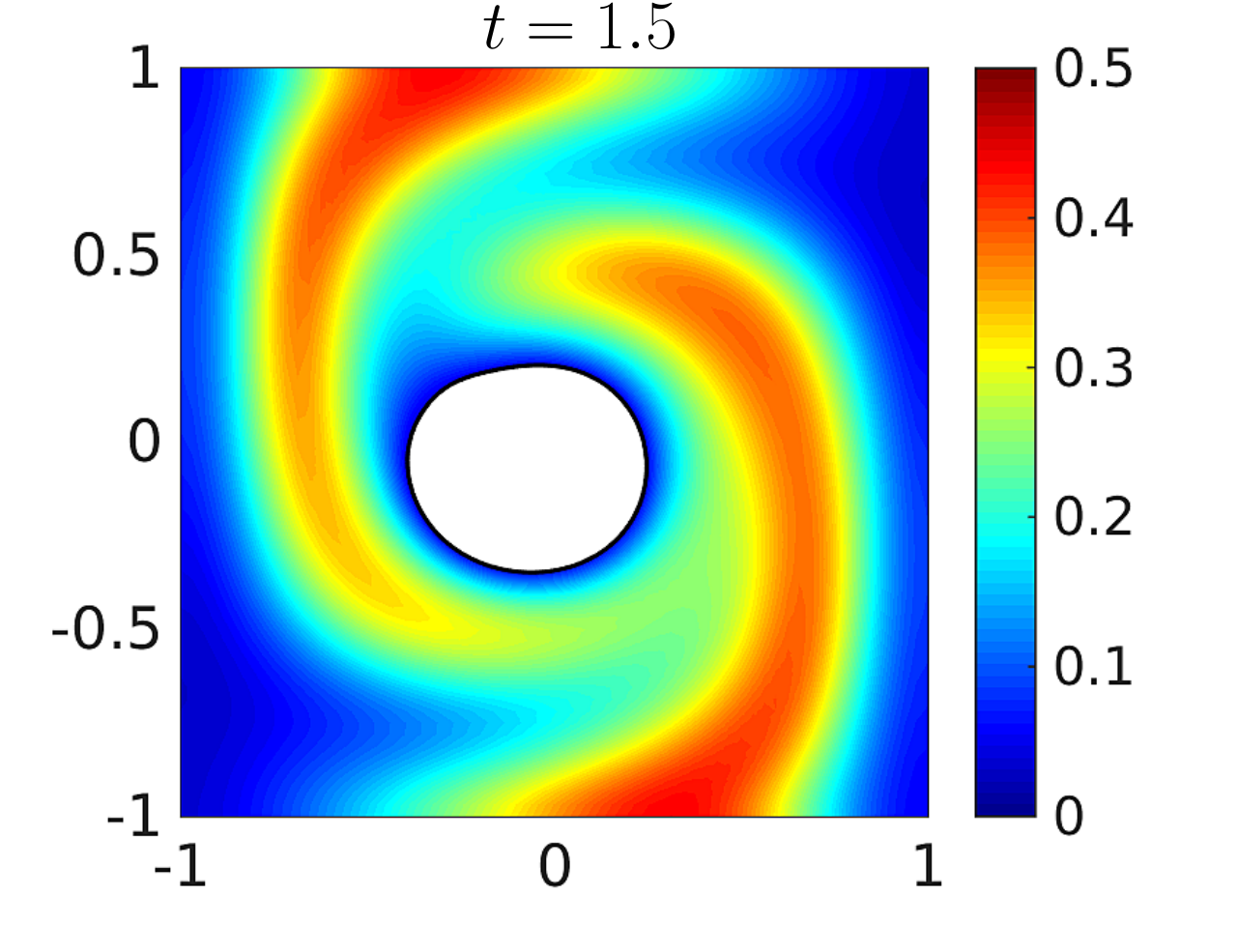}
\includegraphics[width=0.45\textwidth]{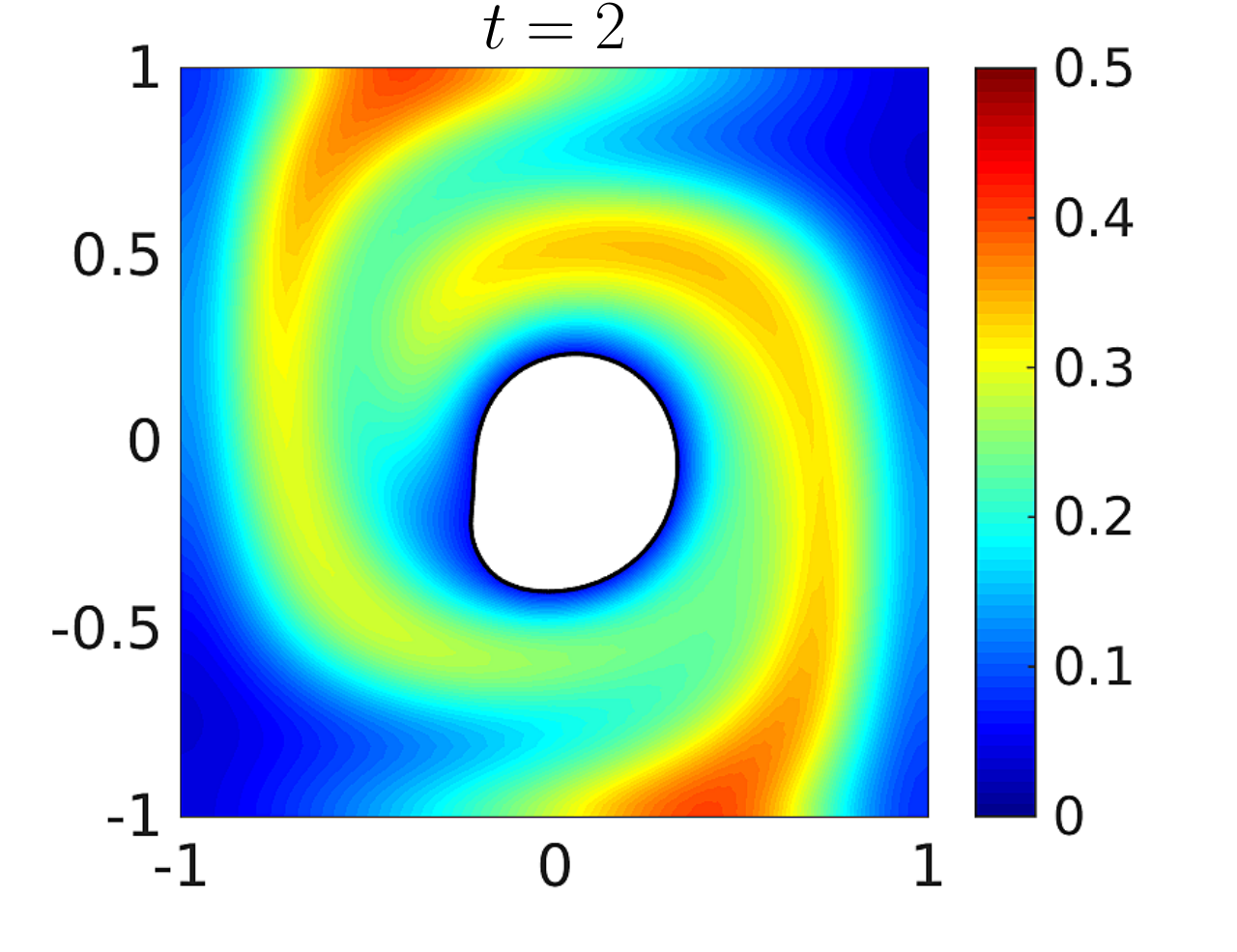}
\caption{Position of the interface and the bulk concentration at time $t=0.5, 1, 1.5, 2$ for mesh size $h=2/64=0.03125$ and time step size $k=h/8$. Results from~\cite{HLZ16ST}. \label{fig:uBLai}}
\end{center}
\end{figure}

\begin{figure}
\begin{center} 
\includegraphics[width=0.45\textwidth]{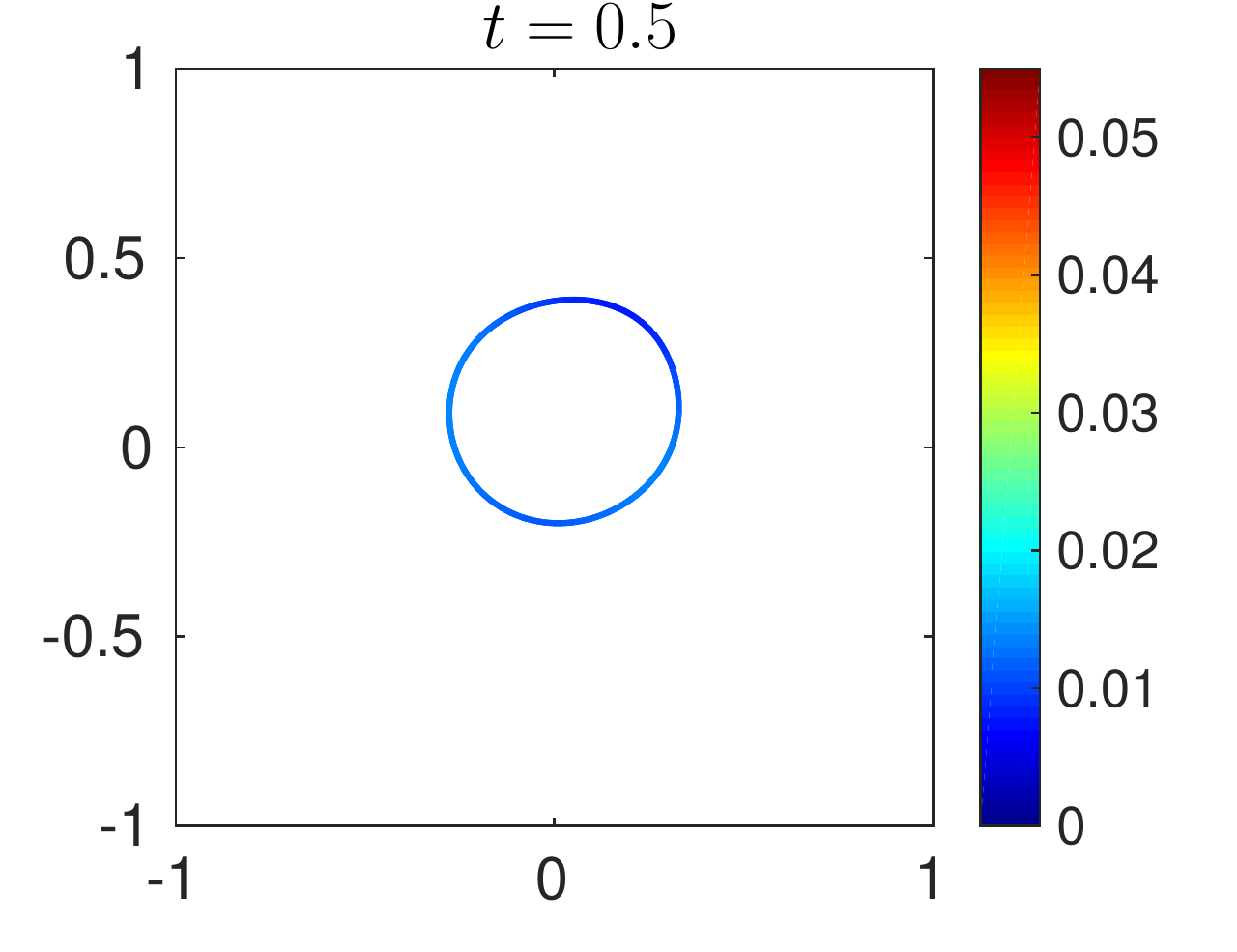}
\includegraphics[width=0.45\textwidth]{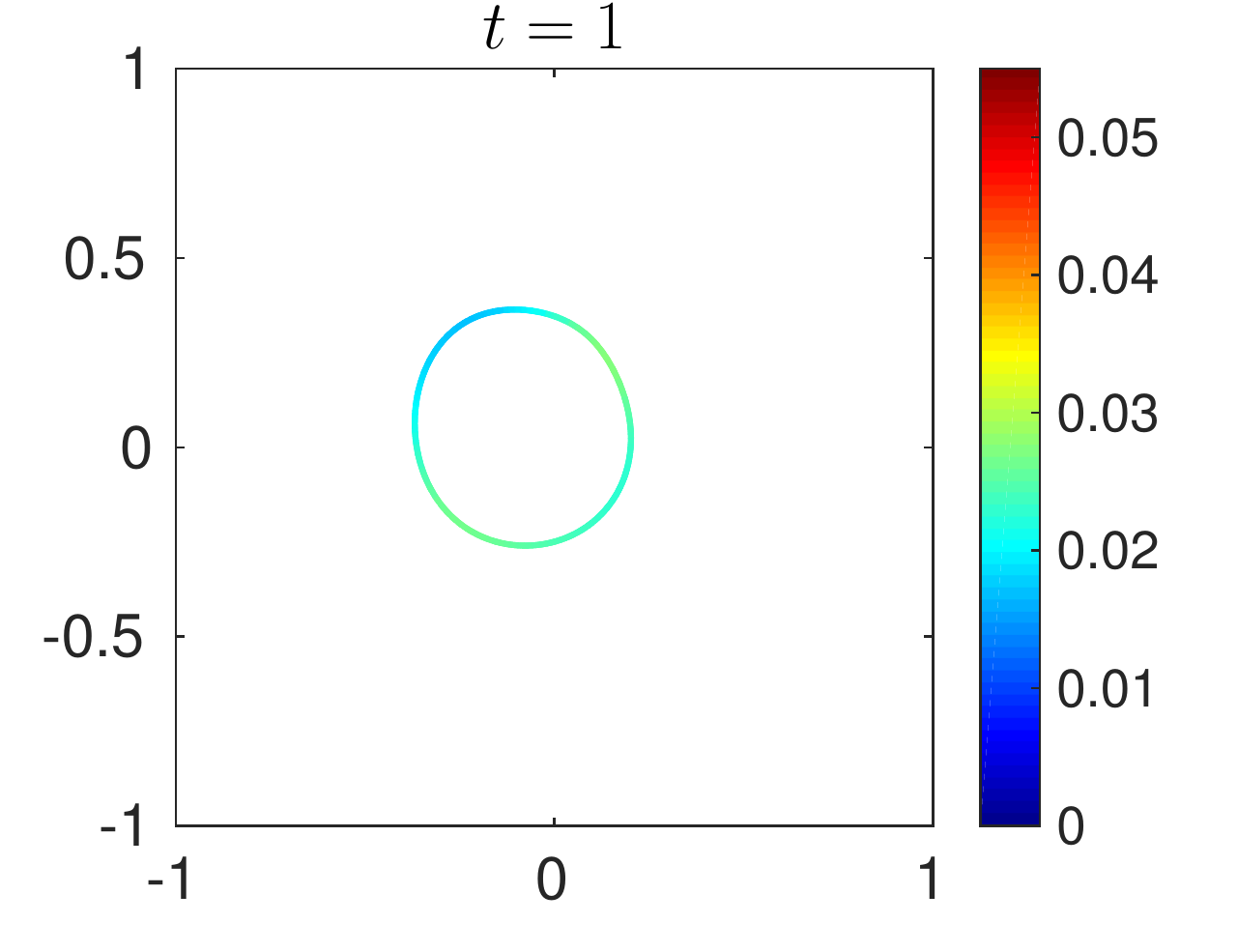}\\
\includegraphics[width=0.45\textwidth]{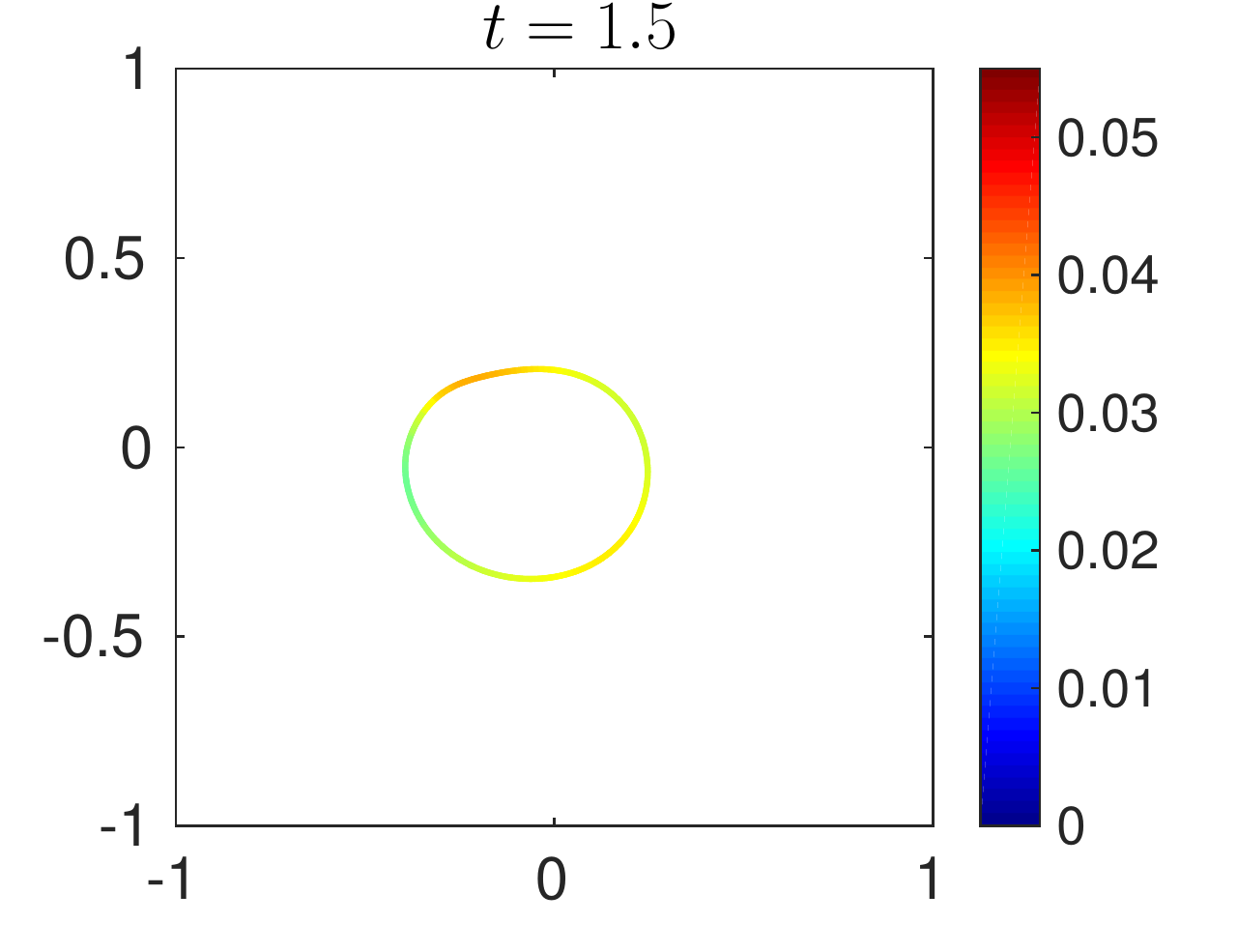}
\includegraphics[width=0.45\textwidth]{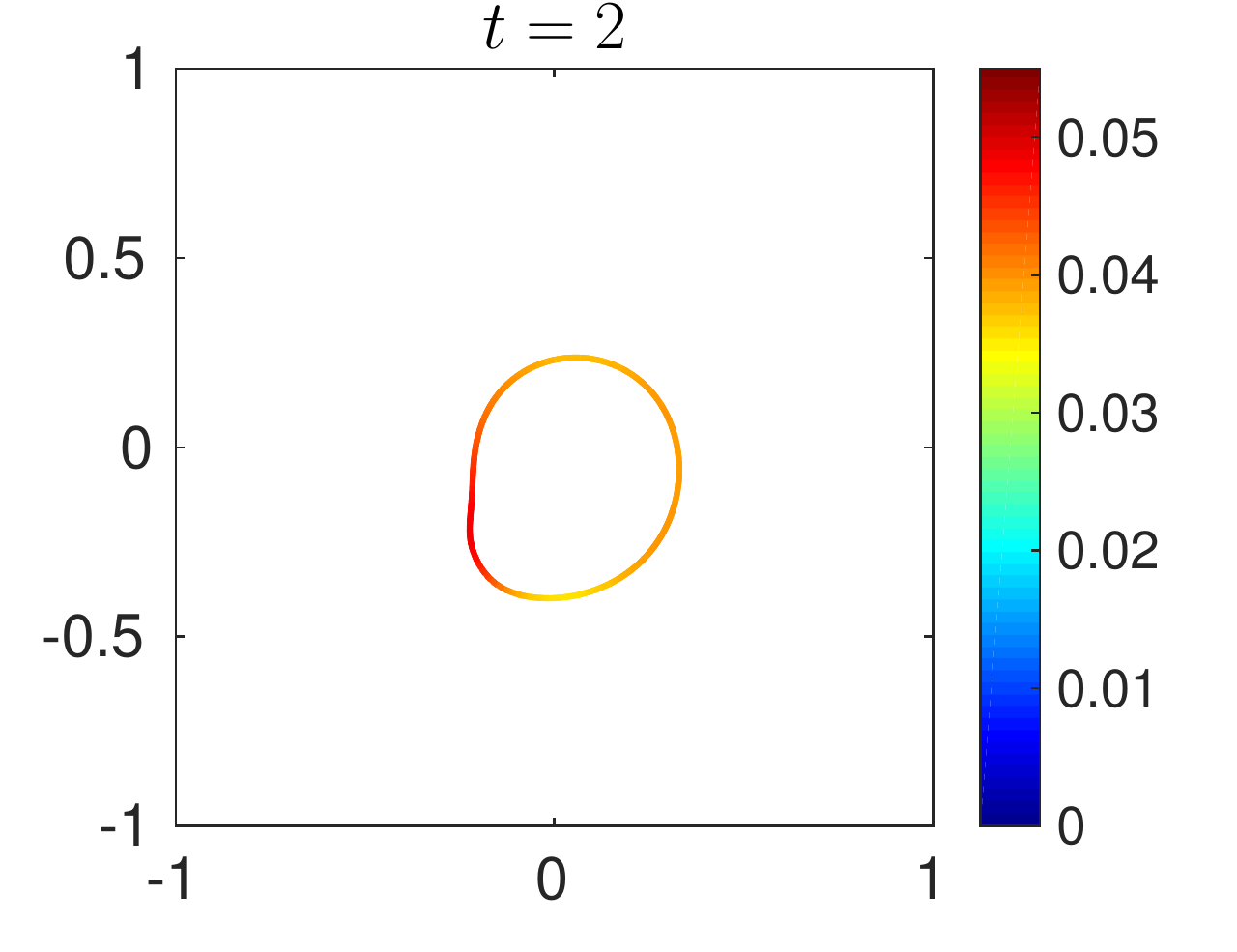}
\caption{Position of the interface and the surface concentration at time $t=0.5, 1, 1.5, 2$ for mesh size $h=2/64=0.03125$ and time step size $k=h/8$.  Results from~\cite{HLZ16ST}.\label{fig:uSLai}}
\end{center}
\end{figure}

\begin{figure}
\begin{center} 
\includegraphics[width=0.5\textwidth]{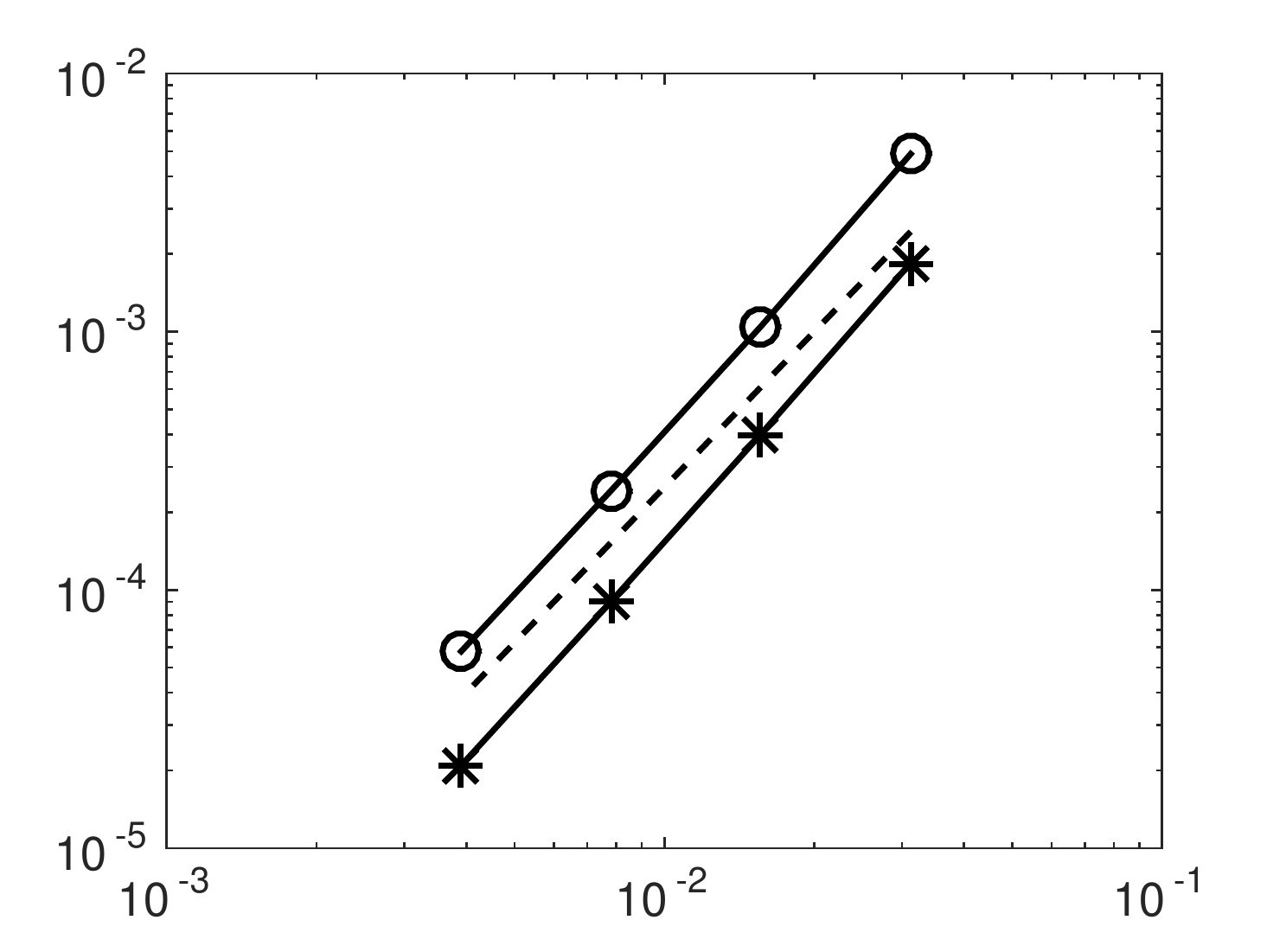} 
\caption{The error $\| (u_{B,h} - u_{B,2h}) \|_{\Omega_{h,1}(0.5)}$ (circles) and $\| (u_{S,h} - u_{S,2h}) \|_{\Gammah(0.5)}$ (stars) measured in the $L^2$ norm versus mesh size $h$.  The dashed line is proportional to $h^2$.  The time step size is $k=0.625h$. Results from~\cite{HLZ16ST}.
 \label{fig:convLai}}
\end{center}
\end{figure}

The bulk and surface surfactant concentrations at times $t=0.5, 1, 1.5, 2$ are shown in Fig.~\ref{fig:uBLai} and Fig.~\ref{fig:uSLai}, respectively. The mesh size is $h=2/64=0.03125$.  We show the error $\| (u_{B,h} - u_{B,2h}) \|_{\Omega_{h,1}(0.5)}$  (represented by circles) and $\| (u_{S,h} - u_{S,2h}) \|_{\Gammah(0.5)}$ (represented by stars) measured in the $L^2$ norm in Fig.~\ref{fig:convLai}. We observe the optimal order of convergence which is second order since we use linear elements in both space and time.
We have measured the order of convergence by using consecutive refinements of the underlying mesh and study $\| (u_{B,h} - u_{B,2h}) \|_{\Omega_{h,1}(0.5)}$  and $\| (u_{S,h} - u_{S,2h}) \|_{\Gammah(0.5)}$. This is also how the convergence is studied in~\cite{ChLai14}.  The method in~\cite{ChLai14} is first order accurate. The errors reported in~\cite{ChLai14} for the bulk concentration, $\| (u_{B,h} - u_{B,2h}) \|_{\Omega_{h,1}(0.5)}$, are smaller for the two coarsest meshes compared to errors using our proposed method but we obtain smaller errors for the two finest meshes. However, for the mesh sizes shown in the figure the errors in the interfacial surfactant concentration, $\| (u_{S,h} - u_{S,2h}) \|_{\Gammah(0.5)}$, reported in~\cite{ChLai14} are smaller than the errors we obtain. This can be understood by the fact that the interface approximation is more accurate in~\cite{ChLai14} where a set of Lagrangian markers are used.

\begin{figure}
\begin{center} 
\includegraphics[width=0.4\textwidth]{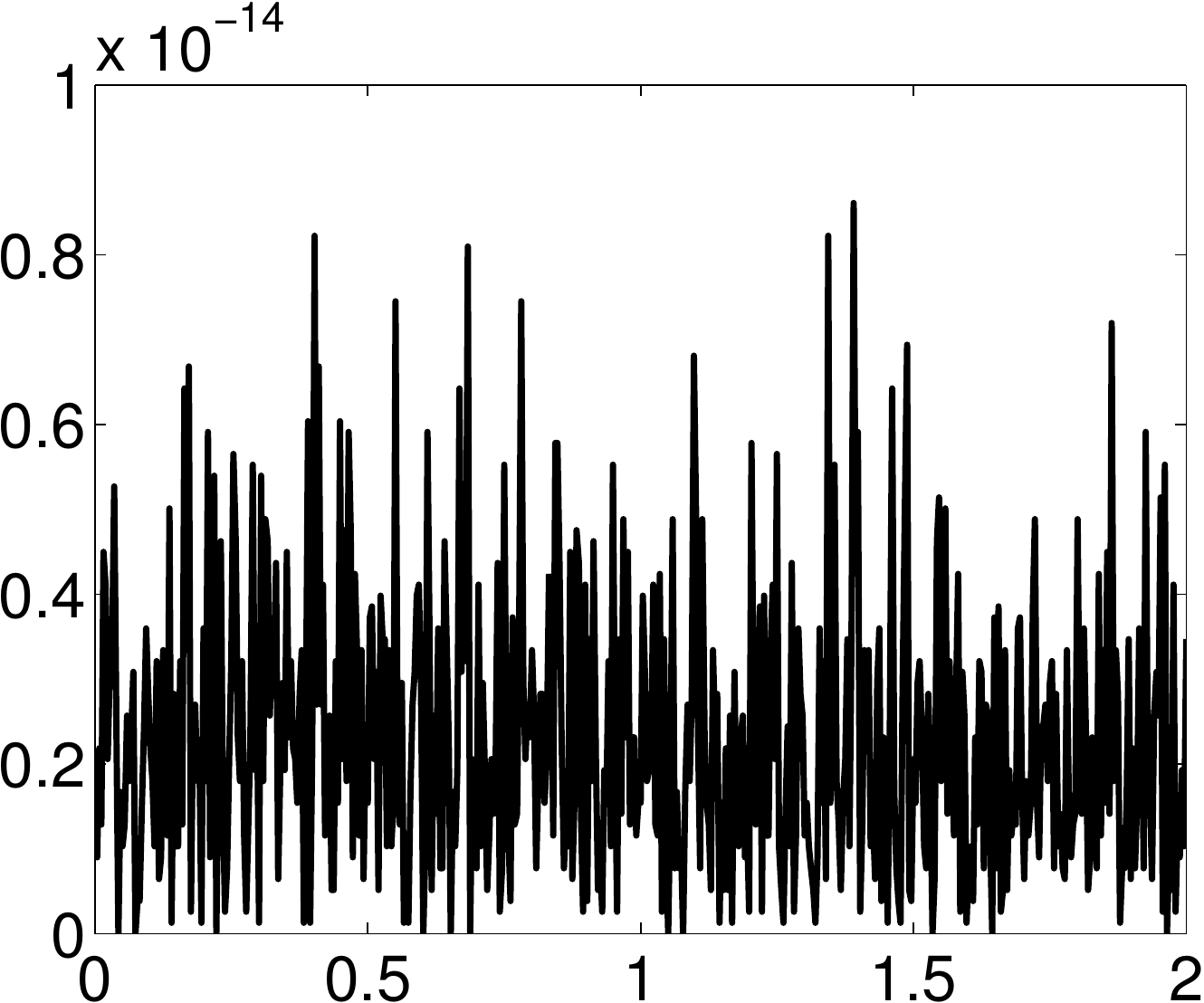} \hspace{0.2cm}
\includegraphics[width=0.4\textwidth]{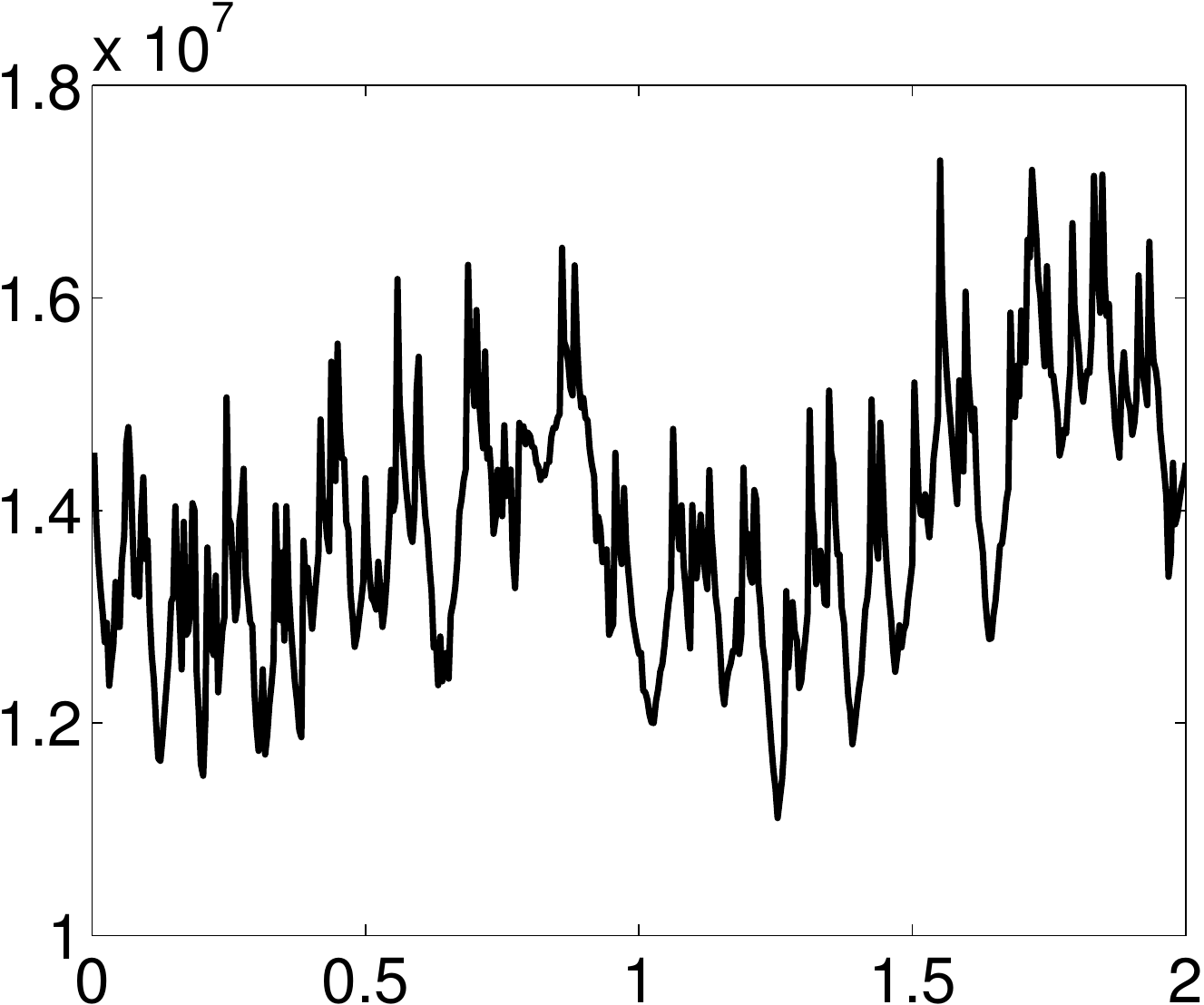}
\caption{Left: The relative error in the total surfactant mass versus time. 
Right: Condition number versus time. 
The mesh size is $h=2/64=0.03125$ in both figures. Results from~\cite{HLZ16ST}. \label{fig:conservandcondLai}}
\end{center}
\end{figure}

In Fig.~\ref{fig:conservandcondLai} we see that the total surfactant mass is conserved. In~\cite{ChLai14} a regularized indicator function is used to extend the bulk equation from $\Omega_1$ to the whole domain. Therefore there is a mass leakage to the domain $\Omega_2$ of the order of the regularization parameter. 
Fig.~\ref{fig:conservandcondLai} also shows the condition number versus time and we see that as the interface evolves the condition number is bounded, independently of how the interface cuts through the mesh.


\subsection{Discussion}\label{sec:discussion}
We studied the space-time CutFEM method developed in~\cite{HLZ16ST} for coupled bulk surface problems modeling the evolution of soluble surfactants. Continuous piecewise linear elements in space and discontinuous piecewise linear elements in time were used and the numerical results show that the method is second order accurate both in space and time. The condition number stays bounded independently of the position of the interface relative to the background mesh due to the face stabilization term that is added in the weak form. The errors we obtain are dominated by the approximation of the interface and we expect to improve the results using a better interface representation. 

A Lagrange multiplier was used to impose the condition~\eqref{eq:conservN} and we therefore had good conservation of the total surfactant mass.  We may consider the same method without the Lagrange multipliers $\lambda$ and $\mu$ as we did for the surface problem. This method is also 
of optimal convergence order but the conservation of the total 
mass of surfactants is lost. Strong imposition of the conservation 
law using Lagrange multipliers essentially compensates for 
numerical errors such as the error in the area of the surface 
and in the volume of the bulk domain, during each time step. 
We also note that using the Reynolds transport theorem one can rewrite the weak form into a conservative form for which condition~\eqref{eq:conservN} is fulfilled at the nodes $t_n$ in the time interval.

To achieve higher order convergence, as for the surface problem, we have to use higher order elements in both space and time in the discretization of the bulk-surface problem (\ref{eq:uBPN}-\ref{eq:uSPN}) and higher order methods for the representation and evolution of the interface. In elements that are cut by the interface we would need quadrature methods for integration in curved domains. Recently, several methods for integration on such curved domains, when the interface is defined implicitly by a level set function, have been proposed, see \cite{Fr15, Le16, Sa15}.  For higher order elements, as we saw in the previous section, other stabilization terms need to be used that can control the condition number.


\end{document}